\documentclass[11pt,leqno]{article}

\usepackage{vmargin}

\usepackage[latin1]{inputenc}
\usepackage[english]{babel}
\usepackage{amsmath}
\usepackage{amssymb}

\usepackage{tabularx}
\usepackage{epsfig, floatflt,amssymb}
\usepackage{moreverb} 
\usepackage{multirow} 
\usepackage{url} 
\usepackage[all]{xy} 
\usepackage{shorttoc} 
\usepackage{textcomp} 
\usepackage{float}
\usepackage[right]{eurosym}

\makeindex 

\usepackage{latexsym}
\usepackage{amssymb,amsbsy,amsmath,amsfonts,amscd}
\usepackage{vmargin}
\usepackage{epsfig}
\usepackage{amsmath}
\usepackage{amsfonts}
\usepackage{amssymb}
\usepackage{graphicx}%
\usepackage[dvips]{color}

\newtheorem{lemme}{Lemma}[section]

\newtheorem{Theorem}{Theorem}[section]

\newtheorem{Remark}{Remark}[section]

\newtheorem{Corollary}{Corollary}[section]

\newcommand{ \vit}{\vec{u}}

\newcommand{\R}{\ifmmode{{\rm I} \hskip -2pt {\rm R}}
    \else{\hbox{$I\hskip -2pt R$}}\fi}

\newcommand{\NN}{\mathbb{N}}

\newcommand{\RR}{\mathbb{R}}

\newcommand{\tore}{\mathbb{T}_3}
\newcommand{\Z}{\mathbb{Z}}
\newcommand{\N}{\mathbb{N}}

\newcommand{\g} {\nabla }

\newcommand{\x} {{\bf x}}

\newcommand{\BEQ} {\begin{equation} }
\newcommand{\EEQ} {\end{equation} }
\makeatletter
\@addtoreset{equation}{section}

\newcommand{\expk}{e^{ {\rm i} \, {\bf k} \cdot \x}}

\usepackage{epsfig}

\newcounter{taskcounter}[section]

\newcounter{technique}[section]

\newcommand{\B}{\mathcal{B}}

\newcommand{\maps}{\rightarrow}

\newcommand{\ip}[2]{\left<#1,#2\right>}

\def\vec#1{\boldsymbol{#1}}

\def\R{\mathcal{R}}

\def\b0{\vec{0}}

\def\bw{\vec{w}}

\def\bu{\vec{u}}

\begin{document}%

\begin{titlepage}
\title{\vskip-2cm\textbf{On a critical Leray$-\alpha$ model of turbulence	}\vskip0.5cm}

\author{\begin{Large}$\hbox{Hani Ali}\thanks{IRMAR, UMR 6625,
Universit\'e Rennes 1,
Campus Beaulieu,
35042 Rennes cedex
FRANCE;
hani.ali@univ-rennes1.fr}$ \end{Large}}
\end{titlepage}
\date{}
\maketitle

\begin{center}
 \textbf{Abstract}
\end{center}

\ \ This paper  aims to study a family of  Leray-$\alpha$ models  with periodic bounbary conditions. These models are  good approximations for the Navier-Stokes equations. We focus our attention on the critical value of regularization ``$\theta$"  that garantees the global well-posedness for these models.  
 We conjecture that $\theta=\frac{1}{4}$ is the  critical value to obtain such   results.
  When alpha goes to zero, we prove that the Leray-$\alpha$ solution, with critical regularization, gives rise to a suitable solution to the Navier-Stokes equations.
 We also introduce an interpolating deconvolution operator that depends on ``$\theta$". Then we extend our results of existence, uniqueness and convergence to a family of regularized magnetohydrodynamics equations. 
 \\
 
MCS Classification : 35Q30, 35Q35, 76F60, 76B03.

\medskip


\section{Introduction}
The dynamics of fluids provide various
highly challenging theoretical, as well as experimental and computational problems
for engineers, physicists and mathematicians. It is widely believed that all
the informations about turbulence are contained in the dynamics of the solutions
of the Navier-Stokes equations (NSE) for viscous, incompressible, homogenous
fluids.\\
The NSE for a homogenous incompressible fluid are usually written as 
\begin{equation}
 \left\{
\begin{array} {llll} \displaystyle
 \frac{\partial \vec{u}^{}}{\partial t}+{\vec{u}^{}} \cdot \nabla
\vec{u}^{} - \nu \Delta \vec{u}^{} + \nabla
p^{} = \vec{f},\\
\nabla \cdot \vec{u}^{}=0,\\
\vec{u}^{}_{t=0}=\vec{u}_{0}^{}.
\end{array}\right.
\end{equation}
The unkowns are  the  velocity vector field $\vec{u}$
and the scalar
pressure  $p$. The viscosity $\nu$, the initial velocity vector field $\vec{u}_{0}$ and  the external
force $\vec{f}$ are given.\\
The dynamics of several conducting incompressible fluids in presence of a magnetic field are described by the  magnetohydrodynamics equations (MHD). The MHD involve coupling Maxwell's equations governing the magnetic field and the NSE governing the fluid motion. The system has the following form
\begin{equation}
 \left\{
\begin{array} {llll} \displaystyle
 \partial_t\bu-\nu_1 \Delta \bu +({\bu}\cdot\nabla)\bu-({\B}\cdot\nabla)\B
+\nabla p+\frac{1}{2}\nabla|\B|^2=0,
\\ \displaystyle
\partial_t\B -\nu_2 \triangle \B +({\bu}\cdot\nabla)\B-({\B}\cdot\nabla)\bu=0,
\\ \displaystyle
\nabla\cdot \B  =\nabla\cdot \bu = 0,\\
\displaystyle
 \B_{t=0}=\B_{0},\; \bu_{t=0}=\bu_{0}.
 \end{array} \right.
\end{equation}
The unkowns are  the  velocity vector field $\vec{u}$, the magnetic vector field $\B$
and the scalar
pressure  $p$. The kinematic viscosity $\nu_1$, the magnetic diffusivity $\nu_2$, the initial velocity vector field $\vec{u}_{0}$  and the initial magnetic vector field $\B_{0}$ are given.\\
 The problem of the global existence
and uniqueness of the solutions of the three-dimensional Navier-Stokes
 equations are among the most challenging problems of contemporary
mathematics. The first attempt was done to Leray \cite{JL34} who established the existence of global weak dissipative solutions  to NSE. Such solutions are called "turbulent" and allow singularities in the velocity fields.  The extension  to no-slip boundary condition was done by Hopf \cite{hopf51}. 
Short time regularity of solutions has also been known for many years, as they have various interesting partial and conditional regularity results \cite{L69,SE63,CF88, RT83}. 
Historically, the partial regularity result  of Scheffer \cite{S76} proving that the potentially singular set in time has zero half-dimensional Hausdorff measure, led to that of Caffarelli, Khon and Nirenberg \cite{CaKoNi82}. 
In \cite{CaKoNi82} they worked with a particular class of weak solutions of the Navier-Stokes equations called  suitable weak solutions. By a suitable weak solution they mean a weak solution of the NSE such that for all $T \in (0,+\infty]$ and for all  non negative fonction  $\phi \in C^{\infty}$  compactly supported in space and time, the following inequality is valid:
\begin{equation}
\label{local 1}
\begin{array}{cccc}
\displaystyle 2\nu\int_{0}^{T}\int_{\tore}|\nabla \vec{u}_{}|^{2}\phi \ d\vec{x}dt \le \displaystyle \int_{0}^{T}\int_{\tore}| \vec{u}_{}|^{2}\left(\phi_t + \nu \Delta \phi \right) + \left(| \vec{u}_{}|^{2} {\vec{u}_{}}+ 2p\vec{u}_{}\right) \cdot \nabla \phi   \  d\vec{x}dt\\
 \hskip 4cm \displaystyle + \displaystyle 2\int_{0}^{T}\int_{\tore} \vec{f} \vec{u}\phi \  d\vec{x}dt 
 \end{array}
\end{equation}
The inequality (\ref{local 1}) is known as local energy inequality. 
Caffarelli, Khon and Nirenberg \cite{CaKoNi82} proved that for the NSE the singular set of a suitable weak solution has parabolic Hausdorff dimension at most equal to $1$. 
This implies that if singularities do exist, they must be relatively rare.\\
Let us mention that it is not known whether  the weak solutions satisfy the local energy inequality. Whereas it is remarkable that the weak solutions constructed by Leray \cite{JL34}, by regularizing the nonlinear term,  or  by adding hyperviscosity are actually suitable \cite{BCI07,DR00,Gu2006}.\\   
 These problems are also open and very important for the MHD equations \cite{DuLi72,RT83mhd} although the MHD equations are not on the Clay Institute's list of prize problems.\\
In order to study various aspects of turbulence, and motivated by the difficulties
related to the complexity of the 3D Navier-Stokes equations, simpler models and
simplifications of hydrodynamics equations have been proposed over the years for exemple \cite{Lions59,lad67, CFHOTW99b,FHT01,FDT02,CHOT05,ILT05,CLT06,OT2007,LL03,LL06,RL}.

\ In this paper,  we study a  Leray-$\alpha$  approximation for the Navier-Stokes equation subject to space-periodic boundary conditions.

The Leray-$\alpha$ equations  we are considering are
\begin{equation}
\label{alpha ns}
 \left\{
\begin{array} {llll} \displaystyle
 \frac{\partial \vec{u}^{}}{\partial t}+\overline{\vec{u}^{}} \cdot \nabla
\vec{u}^{} - \nu \Delta \vec{u}^{} + \nabla
p^{} = \vec{f} \ \ \ \ \ \hbox{in}\ \RR^{+}\times\mathbb{T}_3,\\
\vec{u}=\alpha^{2\theta}(-\Delta)^{{\theta}{}}\overline{\vec{u}}+\overline{\vec{u}},\ \ \hbox{in} \ \mathbb{T}_3  \\
\nabla \cdot \vec{u}^{}= \nabla \cdot \overline{\vit}=0,\
 \displaystyle \int_{\mathbb{T}_3} \vec{u}^{}=\displaystyle \int_{\mathbb{T}_3}\overline{\vit}=0,\\
\vec{u}^{}(t,\vec{x}+{L}\vec{e}_{\vec{j}})=\vec{u}^{}(t,\vec{x}),\\
\vec{u}^{}_{t=0}=\vec{u}_{0}^{}.
\end{array}\right.
\end{equation}
Where $\left\{\vec{e}_{\vec{j}}, \hbox{j=1,2 or 3} \right\}$ is the
canonical basis of $\RR^{3}$, $ L > 0$. \\
We consider these equations on the  three dimensional torus $\mathbb{T}_3=\left ( \RR^3 / {\cal T}_3 \right) $ where ${\cal T}_{3} = 2 \pi \Z^{3} /L$ , $\vec{x} \in \mathbb{T}_3, \ \hbox{and} \ t \in [0,+\infty[.$
The unkowns are  the  velocity vector field $\vec{u}$
and the scalar
pressure  $p$. The viscosity $\nu$, the initial velocity vector field $\vec{u}_{0}$ and the external
force $\vec{f}$ with $\nabla \cdot \vec{f} = 0$ are given.\\
Let $\theta$ be  a non-negative parameter.
The nonlocal operator $(-\Delta)^{{\theta}{}}$ is
defined through the Fourier transform
\begin{equation}
\widehat{(-\Delta)^{{\theta}{}}{\vec{u}({\vec{k}})}}=|{\vec{k}}|^{2\theta}\widehat{\vec{u}}({\vec{k}}).
\end{equation}
Fractionnal order Laplace operator has been used in  another $ \alpha$ model of turbulence in  \cite{OT2007}.
Existence and uniqness of solutions of other modifications of the Navier-Stokes equations have been studied by Ladyzhenskaya \cite{lad67} Lions \cite{Lions59}, M\'{a}lek et al. \cite{mnrr96}.

Motivated by the above work  \cite{OT2007}, we study here the case  $\theta = \frac{1}{4}$ of (\ref{alpha ns}).\\
 We would like to point out that after having finished writing this paper, we were informed that there exists an 
analysis of a general familly of regularized Navier-Stokes and MHD models \cite{HLT10}. In our paper, we note that the Leray-$\alpha$ family considered is a particular case of the general study  in \cite{HLT10} where the results do not recover the critical case  $\theta=\frac{1}{4}.$ The regularization critical values of the other $\alpha$ models studied in  \cite{HLT10} will be reported in a forthcoming paper.\\
Therefore, the initial value problem we
consider in particular is:
\begin{equation}
\label{alpha ns unquart}
 \left\{
\begin{array} {llll} \displaystyle
 \frac{\partial \vec{u}^{}}{\partial t}+\overline{\vec{u}^{}} \cdot \nabla
\vec{u}^{} - \nu \Delta \vec{u}^{} + \nabla
p^{} = \vec{f} \ \ \ \ \ \hbox{in}\ \RR^{+}\times\mathbb{T}_3,\\
\vec{u}=\alpha^{\frac{1}{2}}(-\Delta)^{\frac{1}{4}}\overline{\vec{u}}+\overline{\vec{u}},\ \ \hbox{in} \ \mathbb{T}_3  \\
\nabla \cdot \vec{u}^{}= \nabla \cdot \overline{\vit}=0,\
 \displaystyle \int_{\mathbb{T}_3} \vec{u}^{}=\displaystyle \int_{\mathbb{T}_3}\overline{\vit}=0,\\
\vec{u}^{}_{t=0}=\vec{u}_{0}^{},
\end{array}\right.
\end{equation}
or equivalently (\ref{alpha ns}) with $\theta =  \frac{1}{4} $ and we are working with periodic boundary conditions.
One of the  main results of this paper  is to establish the global well-posedness of the solution to equations (\ref{alpha ns}) in $L^{{2}}(\mathbb{T}_3)^3$ for $\theta=\frac{1}{4}$.
We conjecture that $\theta=\frac{1}{4}$ is the  critical value to obtain the above result   to eqs. (\ref{alpha ns}). 
 Therefore, smooth solutions of eqs. (\ref{alpha ns}) with $ \theta \ge \frac{1}{4}$ do not develop finite-time singularities.\\
 When $\alpha = 0$, eqs. (\ref{alpha ns}) are reduced to the usual Navier-Stokes equations for incompressible fluids.\\
We recall that $H^{\frac{1}{2}}(\mathbb{T}_3)^3$ is a scale-invariant space for Navier-Stokes equations, i.e. if $ \vec{u}(\vec{x},t)$ is a solution to the Navier-Stokes equations, then so is $ \vec{u}_{\lambda}(\vec{x},t):=\lambda\vec{u}(\lambda\vec{x},\lambda^2 t) $ for any $ \lambda >0$.   Kato
 (\cite{Kato1},\cite{Kato2})  shows the
importance for a functionnal   space invariance by  scaling. 
A lot of
works followed which can be summed up  in the
following way: if the initial data are small in an invariant norm
with respect to the scaling then the solution is smooth for all
time.\\
 Another result of this paper 
 is to establish the global well-posedness of the solution to eqs. (\ref{alpha ns}) in $H^{\frac{1}{2}}(\mathbb{T}_3)^3$ for $\theta=\frac{1}{4}$ and  without smallness conditions on the initial data.\\
 We also discuss the relation between the Leray-$\alpha$ equations with $\theta \ge \frac{1}{4}$  and the NSE. The first result of convergence of a $\alpha$  model to the Navier-Stokes equations is proved in \cite{FDT02} where the authors show the existence of a subsequence of weak solutions of the  $\alpha$ model that converge to a Leray-Hopf weak solution of NSE. We improve their result 
by proving an $L^p$ convergence property  in order to show the convergence of the weak solution to a suitable weak solution to the NSE. 
 Moreover, we extend the theory to other models where we consider a deconvolution type regularisation to the NSE and an extension  to the  MHD  equations is also given.
 The deconvolution operator was introduced by  Layton and Lewandowski \cite{LL08}. We will  define the interpolating deconvolution operator in order to regularize the NSE in section \ref{sec:dec}   and the MHD equations in section \ref{sec:MHD}. As the Leray$-\alpha$ regularization we will show an $L^p$ convergence property to the deconvolution operator, in order to deduce that the deconvolution regularization   will  give rise to suitable weak solutions.\\
The critical values of regularization $\theta=\frac{1}{4}$ holds for the MHD equations and the Navier-Stokes equations so it is natural to ask if   the singularities in MHD equations are similar to the singularities in the Navier-Stokes
equations.  A partial answer is given in \cite{Vyalov2008}, where an analog of the known Caffarelli, Khon and Nirenberg result is established to the MHD equations. The extending  of  the partial regularity
result of Scheffer \cite{S76} to the MHD equation is studied in \cite{AA11}.\\

\ The paper is organized as follows. In section 2 we recall some preliminary results  used later in the proof. In section 3, we present the main  existence results (Theorem \ref{2} and Theorem \ref{1} below ). Sections 4 and 5 are devoted to the proof of these two results where we use the standard Galerkin approximation.  In section 6  we show an  $L^p$ convergence result   in order to prove that the Leray$-\alpha$ equations give rise to a suitable solution to the Navier-Stokes equations. In section 7 we introduce the interpolating deconvolution operator and we give  applications to the Navier-Stokes equations.\\
Finally, we study a deconvolution regularization to the MHD equations.

\section{Functional setting and preliminaries}
In this section, we introduce some preliminary material and notations which are commonly used in the mathematical study of fluids, in particular in the study of the Navier-Stokes equations (NSE).  For a more detailed discussion of these topics, we refer to \cite{CF88, DG95, FMRT01, RL09, RT83}.

We denote by $L^p(\tore)^3$ and $H^s(\tore)^3$ the usual Lebesgue and Sobolev spaces over $\tore$.
 For the $L^2$ norm and the inner product, we write $\|\cdot\|_{L^2}$ 
 and $(\cdot,\cdot)_{L^2}$.
  The norm in $H^1(\tore)^3$ is written as $\|\cdot\|_{H^1}$ and its scalar product as $(\cdot,\cdot)_{H^1}$.  We denote by $(H^1(\tore)^{3})'$ the dual space of $H^1(\tore)^3$.
  Note that we have the continuous embeddings
\begin{equation}\label{embed}
H^1(\tore)^3\hookrightarrow  L^2 (\tore)^3 \hookrightarrow (H^{1}(\tore)^3)'.
\end{equation}
Moreover, by the Rellich-Kondrachov compactness theorem (see, e.g., \cite{Evans1998}), these embeddings are compact.\\
Since we work with periodic boundary conditions, we can  characterize the divergence free spaces by using the Fourier series on the 3D torus $\tore$. 
We expand the velocity in Fourier series as:
$$ \displaystyle \vec{u}(t,\vec{x})=\sum_{\vec{k} \in \vec{I}}\hat{\vec{u}}_{\vec{k}}\exp\left\lbrace   i\vec{k}\cdot \vec{x}\right\rbrace ,$$
where $\vec{I}=\left\lbrace \vec{k}  \ \hbox{such  that}  \   \displaystyle \vec{k}=\frac{2\pi \vec{a}}{\vec{L}},\  \vec{a} \in \vec{Z}^{3} , \vec{a} \neq 0 \right\rbrace ,$\\
where for the vanishing space average case, we have the  condition 
$$\hat{\vec{u}^{}}_{0} = \int_{\tore}\vec{u}=0.$$ 
Since the vector $\vec{u}(t,\vec{x})$ is real-valued, we have 
$$ \hat{\vec{u}}_{-\vec{k}}={\widehat{\vec{u}}_{\vec{k}}}^{\ast}\ \ \hbox{for every} \ \vec{k}, $$ 
where $\displaystyle{\widehat{\vec{u}}_{\vec{k}}}^{\ast}$ denotes the complex conjugate of $\hat{\vec{u}^{}}_{\vec{k}}$.\\
In the Fourier space, the divergence-free condition is
$$\hat{\vec{u}^{}}_{\vec{k}}\cdot \vec{k} =0 \ \ \hbox{for all}\  \vec{k}.$$  
For $s \in \RR$, the usual Sobolev spaces $H^s(\tore)^3$  with zero space average can be represented as 
$$\vec{H}^{s}=\{\vec{u}=\sum_{\vec{k} \in \vec{I}}\hat{\vec{u}}_{\vec{k}}\exp{\left\lbrace   i\vec{k}\cdot \vec{x}\right\rbrace }, \ \hat{\vec{u}}_{-\vec{k}}={\widehat{\vec{u}}_{\vec{k}}}^{\ast}, \
\| \vec{u} \|_{\vec{H}^{s}}^{2}< \infty          \},$$
where 
$$\| \vec{u} \|_{\vec{H}^{s}}^{2}= \sum_{\vec{k} \in \vec{I}^{}}|\vec{k}|^{2s} |\hat{\vec{u}}_{\vec{k}}|^{2}.$$

Let 
$$\vec{V}^{s}=\{ \vec{u} \in \vec{H}^{s},\
\hat{\vec{u}^{}}_{\vec{k}}\cdot \vec{k} =0 \},$$

we identify the continuous dual space  of   $\vec{V}^{s'}$ as $\vec{V}^{-s}$ with the pairing 
given by  $$(\vec{u},\vec{v})_{\vec{V}^{s}} = \sum_{\vec{k} \in \vec{I}^{}} |\vec{k}|^{2s}\hat{\vec{u}}_{\vec{k}}\cdot \hat{\vec{v}}_{\vec{-k}} .$$

We denote  $\vec{V}^{0}$ by $ \vec{H} $  and  $\vec{V}^{1}$ by $ \vec{V} $, where   the  norm $ \displaystyle \|\cdot
\|_{{\vec{H}}}^{2}\equiv\|\cdot\|_{L^2}^2$ and the  scalar product as $(\cdot,\cdot)_{\vec{H}}\equiv(\cdot,\cdot)_{L^2}.$  The norm in $\vec{V}$ is written as $\|\cdot
\|_{{\vec{V}}}^{2}\equiv\|\cdot\|_{H^1}^2$ and its scalar product as $(\cdot,\cdot)_{H^1}$.

We denote by $P_\sigma$ the Leray-Helmholtz projection
operator and define the Stokes operator $A:=-P_\sigma\triangle$ with
domain $\mathcal{D}(A):=H^2(\tore)^3\cap \vec{V}$. 
For $\vit \in \mathcal{D}(A)$, we have the norm equivalence $\|A \vit\|_{L^2}\cong\|\vit\|_{\vec{V}^2}$.  Furthermore, in our case it is known that $A=-\triangle$ due to the periodic boundary conditions (see, e.g., \cite{CF88,FMRT01}).\\
It is natural to define the powers $A^{s}$ of the Stokes operator in the periodic case as   
$$A^{s}\vec{u}=\sum_{\vec{k} \in \vec{I}^{}}|\vec{k}|^{2s} \hat{\vec{u}}_{\vec{k}}\exp\left\lbrace   i\vec{k}\cdot \vec{x}\right\rbrace.$$
For $\vit \in \mathcal{D}(A^{s/2})$, we have the norm equivalence $\|A^{s/2} \vit\|_{L^2}\cong\|\vit\|_{\vec{V}^{s}}$.\\
In particular for $ s=1$ we recover the space $ \vec{V}^{}$, hence $\mathcal{D}(A^{1/2})= \vec{V}^{}$, (see, e.g., \cite{FMRT01}).\\

Note also that $\vec{H}^{s+\epsilon}$ is compactly embedded in $\vec{H}^{s}$ (resp. $\vec{V}^{s+\epsilon}$ is compactly embedded in $\vec{V}^{s}$) for any $\epsilon >0$, and we have the following Sobolev embedding Theorem (see \cite{Adams}).
\begin{Theorem}
\label{soblevin}
The space $\vec{H}^{1/2}$ is embedded in $L^{3}(\tore)^3$  and  the space $L^{3/2}(\tore)^3$ is embedded in $\vec{H}^{-{1/2}}$
\end{Theorem}

The following result deals with some interpolation inequalities.
\begin{lemme}\label{interpolation}
Let $T>0, 1\le p_1<p<p_2\le\infty, s_1<s<s_2$  and $\eta\in]0,1[$ such that
$$\frac1p=\frac{\eta}{p_1}+\frac{1-\eta}{p_2}\quad\hbox{and}\quad s=\eta s_1+(1-\eta)s_2.
$$
Let ${\vit}\in \displaystyle \bigcap_{i=1}^2 L^{p_i}([0,T],\vec{H}^{s_i}),$ then ${\vit}\in  L^{p}([0,T],\vec{H}^{s})$ and
$$
\|{\vit}\|_{L^p([0,T], \vec{H}^s)}\le C\|{\vit}\|_{L^{p_1}([0,T], \vec{H}^{s_1})}^{\eta}\|{\vit}\|_{L^{p_2}([0,T], \vec{H}^{s_2})}^{1-\eta}.
$$
The result holds true when we work with the spaces $ L^{p_i}([0,T],\vec{V}^{s_i}).$
\end{lemme}

The regularization effect of the nonlocal operator involved in the relation between   $\overline{\vec{u}}$ and ${\vec{u}}$ is described by the following lemma.
\begin{lemme}
\label{regularite}
Let $\theta\in\RR^+,s\in\RR$ and assume that $\vec{u} \in \vec{H}^{\vec{s}}$.   Then $ \overline{\vec{u}} \in
\vec{H}^{{s}+ 2\theta}$  and 
\begin{equation}
\label{1suralpha2theta}
\|\overline{\vec{u}}  \|_{\vec{H}^{{s}+ 2\theta}}  \le \frac{1}{\alpha^{2\theta}} \|\vec{u} \|_{\vec{H}^{{s}}}.
  \end{equation}
\end{lemme}
\textbf{Proof.} 
When $ \displaystyle \bu=\sum_{{\vec{k}}\in \vec{I}}\hat{\vec{u}}_{\vec{k}}\exp\left\lbrace   i\vec{k}\cdot \vec{x}\right\rbrace $, then 
\BEQ \label {GHGGG}  \overline{\bu} =\sum_{{\vec{k}}\in {\vec{I}}} {\frac{\hat{\bu}_{\vec{k}}}{1 + \alpha^{2 \theta} | {\vec{k}}|^{2\theta}}  }\exp\left\lbrace   i\vec{k}\cdot \vec{x}\right\rbrace.
\EEQ
Formula (\ref{GHGGG}) easily yields  the estimate (\ref{1suralpha2theta}).

\begin{Remark}
For $\theta= \frac{1}{4}$,   $ \overline{\vec{u}} \in
\vec{H}^{{s}+\frac{1}{2} }$  and 
$\displaystyle \|\overline{\vec{u}}\|_{\vec{H}^{{s}+\frac{1}{2}}}  \le {\alpha^{-\frac{1}{2}}} \|\vec{u}\|_{\vec{H}^{s}}.
  $
\end{Remark}
Following the notations of the Navier-Stokes equations, we set \begin{equation}\label{Bdef}
B({\vec{u}},\vec{v}):= P_\sigma({\vec{u}}\cdot \nabla\vec{v})
\end{equation} 
for any 
${\vec{u}},\vec{v}$ in $\vec{V}$. 


We list several important properties of $B$ which can be found for example in \cite{FMRT01}.
\begin{lemme}\label{B:prop}
The operator $B$ defined in \eqref{Bdef} is a bilinear form which can be extended as a continuous map $B:\vec{V}\times \vec{V} \maps \vec{V}'$.\\
   For $\vec{u}$, $\vec{v}$, $\vec{w}\in \vec{V}$,
\begin{equation}\label{B:Alt}
 \ip{B(\vec{u},\vec{v})}{\vec{w}}_{\vec{V}^{-1}}=-\ip{B(\vec{u},\vec{w})}{\vec{v}}_{\vec{V}^{-1}},
\end{equation}
and therefore
\begin{equation}\label{B:zero}
 \ip{B(\vec{u},\vec{v})}{\vec{v}}_{\vec{V}^{-1}}=0.
\end{equation}
\end{lemme}
The following result holds true.
\begin{lemme}\label{B:prop 2.3} The bilinear form $B$ defined in (\ref{Bdef})  satisfies the following:\\
(i) Assume that  $ \vec{u} \in \vec{V}^{\frac{1}{2}}$ ,  $\vec{v}$ and  $\vec{w} \in \vec{V}$,  Then the following inequality holds
 \begin{equation}
\begin{array}{cccccccc}
 |(B({\vec{u}_{}}^{} ,\vec{v}_{}^{} ),\vec{w})| &\le& C \|{\vec{u}}^{} \|_{\vec{V}^{\frac{1}{2}}} \|\vec{v}^{}\|_{\vec{V}} \|\vec{w} \|_{\vec{V}}.
 \end{array}
 \end{equation}
 (ii) $B$ can be extended as a continuous map $B : \vec{V}^{\frac{1}{2}} \times  \vec{V}^{\frac{3}{2}} \rightarrow \vec{V}^{-\frac{1}{2}}$. In particular, for every $\vec{u} \in \vec{V}^{\frac{1}{2}} $, $\vec{v} \in  \vec{V}^{\frac{3}{2}}$  the bilinear form $B$
satisfies the following inequalities:
\begin{equation}
\begin{array}{cccccccc}
  \|B({\vec{u}_{}}^{} ,\vec{v}_{}^{} )\|_{\vec{V}^{-\frac{1}{2}} }
 \le C \|{\vec{u}}^{}_{} \|_{\vec{V}^{\frac{1}{2}}} \|\vec{v}_{}^{}\|_{\vec{V}^{\frac{3}{2}}}.
 \end{array}
 \end{equation}
 \end{lemme}
 
 \textbf{Proof}
 (i) 
 We have 
 \begin{equation}
 \begin{array}{cccccccc}
 |(B({\vec{u}_{}}^{} ,\vec{v}_{}^{} ),\vec{w})|& \le  &\left| \displaystyle  \int_{\tore}{\vec{u}}^{} \otimes \vec{v}_{}^{} : \nabla \vec{w} \right|.
 \end{array}
 \end{equation}
 
 By H$\ddot{o}$lder inequality combined with Sobolev embedding theorem we get
 \begin{equation}
 \begin{array}{cccccccc}
 |(B({\vec{u}_{}}^{} ,\vec{v}_{}^{} ),\vec{w})|
 &\le& \|{\vec{u}}^{} \otimes  \vec{v}_{}^{}\|_{L^2} \|\vec{w}\|_{H^1}\\
 &\le&  \|{\vec{u}}^{} \|_{L^3} \|\vec{v}_{}^{}\|_{L^6 }\|\vec{w}\|_{H^1}   \\
 &\le& C \|{\vec{u}}^{} \|_{\vec{V}^{\frac{1}{2}}} \|\vec{v}_{}^{}\|_{\vec{V}}\|\vec{w}\|_{\vec{V}}.
 \end{array}
 \end{equation}
 (ii)
 We have 
\begin{equation}
\begin{array}{cccccccc}
 \|B({\vec{u}_{}}^{} ,\vec{v}_{}^{} )\|_{\vec{V}^{-\frac{1}{2}} }
 &\le& c \|{\vec{u}_{}}^{} \|_{L^{3}} \|\nabla\vec{v}_{}^{}\|_{L ^{3} }\\
 &\le& c \|{\vec{u}}^{}_{} \|_{\vec{V}^{\frac{1}{2}}} \|\vec{v}_{}^{}\|_{\vec{V}^{\frac{3}{2}}}.
 \end{array}
 \end{equation}
Where he have used H$\ddot{o}$lder inequality and the Sobolev injection given in Theorem \ref{soblevin}.\\ 

The following result is given for any $\vec{u} \in L^{\infty}([0,T];\vec{H}^{}) \cap L^{2}([0,T];\vec{V}^{})$.
\begin{lemme}
\label{lemmetermenonlineaire}
If $\vec{u} \in L^{\infty}([0,T];\vec{H}^{}) \cap L^{2}([0,T];\vec{V}^{})$, then $ \overline{\vec{u}} \in L^{\infty}([0,T];\vec{V}^{\frac{1}{2}}) \cap L^{2}([0,T];\vec{V}^{\frac{3}{2}})$ and 
$B(\overline{\vec{u}},\vec{u}) \in L^{2}([0,T];\vec{V}^{-1})$.
\end{lemme} 
\textbf{Proof.}
The regularity on $\overline{\vec{u}}$ is obtained from  Lemma \ref{regularite}.\\
Then we use Lemma \ref{B:prop 2.3} to obtain 
 $$\|(B(\overline{\vec{u}_{}}^{} ,\vec{u}_{}^{} ))\|_{\vec{V}^{-1}}
 \le  C\|\overline{\vec{u}}^{} \|_{\vec{V}^{\frac{1}{2}}} \|\vec{u}_{}^{}\|_{\vec{V}}$$
Now,
$$
\|\overline{\vec{u}}^{} \|_{\vec{V}^{\frac{1}{2}}}  \in L^{\infty}(0,T) \ \ and \ \  \| \vec{u}\|_{\vec{V}} \in L^{2}(0,T),$$ 
implies that
$\displaystyle \|\overline{\vec{u}}^{}\|_{\vec{V}^{\frac{1}{2}}} \|\vec{u}_{}^{}\|_{\vec{V}^{}} \in L^2(0,T) $ (that is $L^{\infty}([0,T];\vec{V}^{\frac{1}{2}}) \cap  L^{2}([0,T];\vec{V}^{}) \subseteq  L^{2}([0,T];\vec{H}^{})$ ).\\



 

We remove the pressure from further consideration by projecting (\ref{alpha ns}) by $P_\sigma$ and searching for solutions in the space $\vec{V}^{s}$.\\
 Thus, we obtain the infinite dimensional evolution equation 
 
 \begin{equation}
\label{alpha ns prime}
 \left\{
\begin{array} {llll} \displaystyle
 \frac{d \vec{u}^{}}{d t}  +\nu A \vec{u}^{} + B(\overline{\vec{u}},\vec{u})
 = \vec{f},\\
\vec{u}=\alpha^{2\theta}A^{{\theta}{}}\overline{\vec{u}}+\overline{\vec{u}},\\
\vec{u}^{}_{t=0}=\vec{u}_{0}^{}.
\end{array}\right.
\end{equation}

 Note that the pressure may be reconstructed from $\overline{\vec{u}}$ and $\vec{u}$ by solving the elliptic 
equation $$\Delta p= \nabla \cdot (\overline{\vec{u}} \cdot \nabla \vec{u}).$$\\
Later in Sect.4  and Sect. 5, we will prove that eqs. (\ref{alpha ns}) have unique regular solutions $(\vec{u}_{},p_{})$
such that the velocity part $\vec{u}_{} $ solve the eqs. (\ref{alpha ns prime}).\\

We conclude this section with the interpolation lemma of Lions and Magenes \cite{JL69} and a compactness Theorem see \cite{sim87}.
 \begin{lemme}
 \label{lionslemma}
 Let $s >0$ and suppose that 
 $$  \vec{u} \in L^{2}([0,T];\vec{V}^{s}) \ \hbox{ and } \  \displaystyle \frac{d\vec{u}^{}}{dt} \in L^{2}([0,T],\vec{V}^{-s}).$$
  Then \\
  (i)  $\vec{u} \in C([0,T];\vec{H})$, with
  $$\sup_{t\in [0,T]} \|\vec{u}(t)\|_{\vec{H}}  \le c \left( \|\vec{u}\|_{L^{2}([0,T];\vec{V}^{s})}  +    \|\frac{d\vec{u}^{}}{dt}\|_{L^{2}([0,T],\vec{V}^{-s})}      \right) $$
  and\\
  (ii) $ \displaystyle \frac{d}{dt}\|\vec{u}\|_{\vec{H}}^{2}= 2\left\langle\frac{d\vec{u}^{}}{dt},\vec{u} \right\rangle_{\vec{V}^{-s},\vec{V}^{s}}$.
  \end{lemme}
 To set more regularity for the unique  solution, we  also   need the following similar result.
 \begin{Corollary}
 Assume that for  some $ k \ge 0$ and  $s >0$,  
 $$  \vec{u} \in L^{2}([0,T];\vec{V}^{k+s}) \ \hbox{ and } \  \displaystyle \frac{d\vec{u}^{}}{dt} \in L^{2}([0,T],\vec{V}^{k-s}).$$
  Then \\
  (i)  $\vec{u} \in C([0,T];\vec{V}^{k})$,  and\\
  (ii) $ \displaystyle \frac{d}{dt}\|A^{k/2}\vec{u}\|_{\vec{H}}^{2}= 2\left\langle A^{k/2}\frac{d\vec{u}^{}}{dt},A^{k/2}\vec{u} \right\rangle_{\vec{V}^{-s},\vec{V}^{s}}$.
 \end{Corollary}
 The result holds true in the space $\vec{H}^{s}.$
 \begin{Theorem}
 \label{compact}
 Let $X \subset\subset H \subset Y$ be Banach spaces, where $X$ is reflexive. Suppose that $\vec{u}_{n}$ is a sequence such that it is uniformly bounded in $L^{2}([0,T];X),$ and 
 $\displaystyle \frac{d\vec{u}_{n}}{dt}$ is uniformly bounded in $L^{2}([0,T],Y)$. Then there is a subsequence which converges strongly in $L^{2}([0,T],H)$.
 \end{Theorem}
 The above Theorem is a special case of the following more general result:
  \begin{Theorem}(Aubin-Lions)
 \label{Aubin-Lions}
 Let $X \subset\subset H \subset Y$ be Banach spaces, where $X$ is reflexive.
Suppose that $\vec{u}_{n}$ is a sequence such that it is uniformly bounded in $L^{p}([0,T];X),$ ($\displaystyle 1<p<\infty$), and 
 $\displaystyle \frac{d\vec{u}_{n}}{dt}$ is uniformly bounded in $L^{q}([0,T],Y)$, ($  \displaystyle  1 \le q \le \infty$). Then $\vec{u}_{n}$ is relatively compact in $L^{p}([0,T],H)$.
 \end{Theorem}

\section{Main existence theorems}

\ In the following, we will assume that $ \alpha >0$, $T>0$. 


One of the aims  of this paper is to establish the global well-posedness of the solution to eqs. (\ref{alpha ns}) in $L^{{2}}(\mathbb{T}_3)^3$ for $\theta=\frac{1}{4}$. It must be mentioned that the result holds true for $\theta \ge \frac{1}{4}$.
We conjecture that $\theta=\frac{1}{4}$ is the  critical value to obtain the following result to eqs. (\ref{alpha ns}). 
\begin{Theorem}
\label{2}
  For any $T>0$,  let $\vec{f} \in L^{2}([0,T],\vec{V}^{-{1}})$ and $\vec{u}_{0} \in \vec{H}$. Assume that  $\theta=\frac{1}{4}$. Then there exists a unique     solution $(\vec{u}_{},p_{}):=(\vec{u}_{\alpha},p_{\alpha})$ to eqs. (\ref{alpha ns}) that satisfies 
$ \vec{u}_{} \in  C([0,T];\vec{H}) \cap L^{2}([0,T];\vec{V})$ and $\displaystyle \frac{d\vec{u}}{dt} \in L^{2}([0,T];\vec{V^{-1}})$ and  ${p}_{} \in  L^{2}([0,T],L^{{2}}(\mathbb{T}_3))$.\\
Such that $\vec{u}_{}$ verifies 
$$\left\langle  \frac{d \vec{u}_{}}{d t}  +\nu A \vec{u}_{} + B(\overline{\vec{u}},\vec{u})
 -\vec{f}, \phi  \right\rangle_{\vec{V}^{-{{1}{}}},\vec{V}^{{{1}{}}}}=0 $$
 for every $\phi \in \vec{V}^{}$ and almost every $ t \in (0,T)$.
Moreover, this solution depends continuously on the initial data $\vec{u}_{0}$ in the $L^{{2}}$ norm.
\end{Theorem}
 Therefore, smooth solutions of eqs. (\ref{alpha ns}) with $ \theta \ge \frac{1}{4}$ do not develop finite-time singularities.\\

We also prove 
  the global well-posedness of the solution to eqs. (\ref{alpha ns}) in $H^{\frac{1}{2}}(\mathbb{T}_3)^3$ for $\theta=\frac{1}{4}$ and  without smallness conditions on the initial data. 
\begin{Theorem}
\label{1}
For any $T>0$, let $\vec{f} \in L^{2}([0,T],\vec{V}^{-{\frac{1}{2}}})$ and $\vec{u}_{0} \in \vec{V}^{\frac{1}{2}}$. Assume that $\theta= \frac{1}{4}$. Then  there exists a unique solution  $(\vec{u}_{},p_{}):=(\vec{u}_{\alpha},p_{\alpha})$ to eqs. (\ref{alpha ns})  that satisfies 
 $ \vec{u}_{} \in  C([0,T];\vec{V}^{\frac{1}{2}}) \cap L^{2}([0,T];\vec{V}^{\frac{3}{2}})$, $\displaystyle \frac{d\vec{u}_{}^{}}{dt} \in L^{2}([0,T],\vec{V}^{{-\frac{1}{2}}})$ and  ${p}_{} \in  L^{2}([0,T],H^{\frac{1}{2}}(\mathbb{T}_3))$.\\
Such that  $\vec{u}$ satisfies 
$$\left\langle  \frac{d \vec{u}^{}}{d t}  +\nu A \vec{u}^{} + B(\overline{\vec{u}},\vec{u})
 -\vec{f}, \phi  \right\rangle_{\vec{V}^{-\frac{1}{2}},\vec{V}^{\frac{1}{2}}}=0 $$
 for every $\phi \in \vec{V}^{\frac{1}{2}}$ and almost every $ t \in (0,T)$.
Moreover, this solution depends continuously on the initial data $\vec{u}_{0}$ in the $\vec{V}^{\frac{1}{2}}$ norm.
\end{Theorem}

We note that the above result holds true for any $\displaystyle s >0 $ and without smallness conditions on the initial data.

 
\section{Proof of Theorem \ref{2}}

The proof is divided into  four steps 
\subsection{Galerkin approximation}
Let us define $$H_{m} \equiv span\left\{ {\exp\left\lbrace   i\vec{k}\cdot \vec{x}\right\rbrace : |\vec{k}| \le m }\right\}.$$ 
We look at the finite-dimensional equation obtained by keeping only the first $m$ Fourrier modes. 
 In order to use classical tools for systems of ordinary differential equations we need that $\vec{f}$ belongs to $C([0,T],\vec{V}^{-1})$.  To do so, we extend $\displaystyle\vec{f}$  outside $[0,T]$ by zero and we set $ \displaystyle \vec{f}_\epsilon=\rho_\epsilon\ast\vec{f}$  where $\displaystyle \rho_\epsilon(t)=\frac{1}{\epsilon}\rho(\frac{t}{\epsilon})$,  $0 \le \displaystyle \rho(s) \le 1, \displaystyle \rho(s)=0$ for $|s| \ge 1$, and $\displaystyle \int_{\RR}\rho=1$. 
So the approximate sequence $\displaystyle \vec{f}_\epsilon$  is very smooth with respect to time  for $\displaystyle \rho$ smooth and converges to $\displaystyle \vec{f}$ in the sense that 
\begin{align}
 {\vec{f}}_\epsilon &\rightarrow  {\vec{f}} \textrm{  strongly in  }
L^2([0,T];\textbf{V}^{-1}) \textrm{ when } \epsilon \rightarrow 0. \label{conv1}
\end{align} 

The Galerkin approximation of eqs. (\ref{alpha ns prime}) with $\theta = \frac{1}{4}$ is given by 

\begin{equation}
\label{galerkine}
 \left\{
\begin{array} {llll} \displaystyle
 \frac{d \vec{u}^{m}_{}}{d t}  +\nu P_m A \vec{u}^{m}_{} +P_{m} B(P_m\overline{\vec{u}_{}}^{m},P_m\vec{u}^{m}_{})
 = P_{m}\vec{f}_{\frac{1}{m}},\\
\vec{u}^{m}_{}=\alpha^{\frac{1}{2}}A^{{\frac{1}{4}}{}}\overline{\vec{u}_{}}^{m}+\overline{\vec{u}}^{m}_{},\\
\vec{u}^{m}_{}(0)=P_{m}\vec{u}_{0}.
\end{array}\right.
\end{equation}

Where for some $m \in \NN$ and all $-1 \le s \le 2$, $P_{m}(\textsl{w})\equiv \sum_{|\vec{k} | \le m}\hat{\textsl{w}}_{\vec{k}}\exp{\left\lbrace   i\vec{k}\cdot \vec{x}\right\rbrace } : \vec{V}^s \rightarrow H_m $ is 
the orthogonal projector  onto the first $m$ Fourier modes that  verifies (see in \cite{mnrr96} for more details):
 \begin{align}
 \| P_m\|_{ \mathcal{L} ( \vec{V}^{s}, \vec{V}^{s})} \le 1, \quad \textrm{and for all } \vec{v} \in \vec{V}^s:  P_m \vec{v} {\rightarrow}  {\vec{v}} \textrm{ strongly in  } \vec{V}^s  \textrm{ when } m \rightarrow \infty.
  \label{pm}
\end{align} 
  
The classical theory of ordinary differential equations implies that eqs. (\ref{galerkine}) have a unique $C^1$ solution $\vec{u}^m $  for a given time  interval that, \textit{a priori}, depends on $m$, such that $\vec{u}^m =P_m \vec{u}$ and $\nabla \cdot \vec{u}^m=0$.  Our goal is to show that the solution remains finite for all positive times, wich implies that $T_m = \infty$.\\
In the next subsection, we find uniform energy estimates for this  solution with respect to  $m$.

\subsection{Energy estimates in $\vec{H}$ }
We  follow here a similar method to the one used  for the Navier-Stokes equations (see \cite{RT83}).
\begin{lemme}
\label{lemme2}
Let $\vec{f} \in L^{2}([0,T],\vec{V}^{{-1}})$ and $\vec{u}_{0} \in \vec{H} $, there exists $K_{1}(T)$ and $K_{2}(T)$ independent of m such that the solution $\vec{u}_{}^{m}$ to the Galerkin truncation \ref{galerkine} satisfies 
\begin{equation}
\label{estimatel2h1}
\|\vec{u}_{}^{m}\|_{L^{2}([0,T],\vec{V})}^2 \le K_{1}(T), \ \hbox{for all} \  T \ge 0   , \ \ \  \hbox{where} \ K_{1}(T)= \frac{1}{\nu}\left( \|\vec{u}_{0}\|^2_{\vec{H}}+\frac{1}{\nu} \|\vec{f} \|_{L^{2}([0,T],\vec{V}^{{-1}})}^{2} \right)      
\end{equation}
and
\begin{equation}
\|\vec{u}_{}^{m}\|_{L^{\infty}([0,T],\vec{H})}^2 \le K_{2}(T), \ \hbox{for all} \  T \ge 0   , \ \ \  \hbox{where} \  K_{2}(T)= {\nu}(   K_{1}(T)   ).
\end{equation}
\end{lemme}   

\textbf{Proof}
Taking the $L^2$-inner product of the first equation of (\ref{galerkine})  with $\vec{u}$ and integrating by parts, using (\ref{B:zero}), the incompressibility of the velocity field and the duality relation we obtain
\begin{eqnarray}
\frac12\frac{d}{dt}\|\vec{u}_{}^{m}\|_{\vec{H}}^2+\nu\|\nabla \vec{u}_{}^{m}\|_{\vec{H}}^2&=&\int_{\tore}P_{m}\vec{f}_{\frac{1}{m}} \cdot \vec{u}_{}^{m}dx
\le \|\vec{f}\|_{\vec{V}^{-1}} \|\vec{u}_{}^{m}\|_{\vec{V}^1}.
\end{eqnarray}
Using Young inequality we get
\begin{eqnarray}
\label{page285temam}
\frac{d}{dt}\|\vec{u}_{}^{m}\|_{\vec{H}}^2+\nu\|\nabla \vec{u}_{}^{m}\|_{\vec{H}}^2&\le&\frac{1}{\nu}\|\vec{f}\|_{{\vec{V}^{-1}}}^2.
\end{eqnarray}
In particular, this  leads to  the estimate
\begin{equation}
\label{contradiction3.6}
\sup_{t \in [0,T_m)}\|\vec{u}_{}^{m}\|_{\vec{H}}^2 \le \|\vec{u}_{0}\|_{\vec{H}}^2+\frac{1}{\nu} \|\vec{f} \|_{L^{2}([0,T],\vec{V}^{{-1}})}^{2}
\end{equation}
This implies that $T_m = T$. Indeed, consider $[0,T_{m}^{max})$ the maximal interval of existence. Either $T_{m}^{max}= T$ and we are done, or $T_{m}^{max} < T$ and we have $\limsup_{t\rightarrow(T_{m}^{max})^{-}} \|\vec{u}_{}^{m}(t)\|_{\vec{H}}^2=\infty$ a contradiction to (\ref{contradiction3.6}). Hence we have global existence of $\vec{u}_{}^{m}$ for  $\vec{f} \in L^{2}([0,\infty),\vec{V}^{{-1}})$ and $\vec{u}_{0} \in \vec{H} $ and hereafter we take an arbitrary interval $[0,T]$ and we assume that $ \vec{f} \in L^{2}([0,T],\vec{V}^{{-1}})$.\\
Integrating (\ref{page285temam}) with respect to  time gives the desired estimate (\ref{estimatel2h1}) for all $t \in [0,T]$.\\ 

Now we use the regularization effect described in Lemma \ref{regularite} to  get the following estimates on $\overline{\vec{u}_{}}^m$.
\begin{lemme}
\label{l3}
Under the notation of Lemma \ref{lemme2},
we have  
\begin{equation}
\|\overline{\vec{u}_{}}^{m}\|_{L^{2}([0,T],\vec{V}^{\frac{3}{2}})}^2 \le \left(\frac{1}{\alpha^{}}\right) K_{1}(T), \ \hbox{for all} \  T \ge 0, 
\end{equation}
and
\begin{equation}
\|\overline{\vec{u}_{}}^{m}\|_{L^{\infty}([0,T],\vec{V}^{\frac{1}{2}})}^2 \le \left(\frac{1}{\alpha^{}}\right) K_{2}(T), \ \hbox{for all} \ T \ge 0.
\end{equation}
\end{lemme}
 By interpolation  (see Lemma \ref{interpolation}) between $L^{2}([0,T],\vec{V}^{\frac{3}{2}})$ and $L^{\infty}([0,T],\vec{V}^{\frac{1}{2}})$, we can deduce the following result
\begin{Corollary} Under the notation of Lemma \ref{lemme2},
we deduce 
\begin{equation}
\|\overline{\vec{u}_{}}^{m}\|_{L^{4}([0,T],\vec{V})}^2\le \left(\frac{1}{\alpha^{}}\right) K_{1}(T)^{1/2} K_{2}(T)^{1/2}, \ \hbox{for  all} \ T \ge 0.
\end{equation}
\end{Corollary}

Our next result provides an estimate on the time derivative of $\vec{u}_{}^{m}$. 
\begin{lemme}
\label{l4}
Let $\vec{f} \in L^{2}([0,T],\vec{V}^{-1})$ and $\vec{u}_{0} \in \vec{H}$, there exists $K_{3}(T)$ independent of m such that the time derivative of  the solution $\vec{u}_{}^{m}$ to the Galerkin truncation (\ref{galerkine}) satisfies 
\begin{equation}
\label{du sur dt 1}
\|\frac{d\vec{u}_{}^{m}}{dt}\|_{L^{2}([0,T],\vec{V}^{-1})}^2 \le K_{3}(T), \ \ \ 
\end{equation}
{where}
$$
  K_{3}(T)= 4\left( \nu^{2} + C \left(1/\alpha^{}\right) K_{2}(T) \right)^{}
K_{1}(T)+ 2\| \vec{f}\|_{L^{2}([0,T],\vec{V}^{-1})}^{2}.
$$
\end{lemme}
\textbf{Proof.}
Taking the $\vec{H}^{-1}$ norm of (\ref{galerkine}), one obtains :  
\begin{equation}
\|\frac{d\vec{u}_{}^{m}}{dt}\|_{\vec{V}^{-1}}  \le  \nu
\|{\vec{u}_{}^{m}}{}\|_{\vec{V}^1}+\| B(\overline{\vec{u}_{}}^{m} ,\vec{u}_{}^{m} )\|_{\vec{V}^{-1}}+ \| \vec{f}\|_{\vec{V}^{-1}}
\end{equation}
where we note that 
$$\| B(\overline{\vec{u}_{}}^{m} ,\vec{u}_{}^{m} )\|_{\vec{V}^{-1}}= \sup{\left\{(B(\overline{\vec{u}_{}}^{m} ,\vec{u}_{}^{m} ),w);  w \in \vec{V}^1, \|w\|_{\vec{V}^1} \le 1 \right\}}.$$
It remains to show that $\| B(\overline{\vec{u}_{}}^{m} ,\vec{u}_{}^{m} )\|_{\vec{V}^{-1}}$  is bounded.
Indeed, we have  from Lemma \ref{B:prop 2.3} that 
\begin{equation}
\begin{array}{cccccccc}
\| B(\overline{\vec{u}_{}}^{m} ,\vec{u}_{}^{m} )\|_{\vec{V}^{-1}}
 &\le& C \|\overline{\vec{u}_{}}^{m} \|_{\vec{V}^{\frac{1}{2}}} \|\vec{u}_{}^{m}\|_{\vec{V}}.
 \end{array}
 \end{equation}

  It follows that 
 \begin{equation}
\|\frac{d\vec{u}_{}^{m}}{dt}\|_{\vec{V}^{-1}}^{2}  \le  4( \nu^{2} + C \|\overline{\vec{u}_{}}^{m} \|_{\vec{V}^{\frac{1}{2}}}^{2})^{}
\|{\vec{u}_{}^{m}}{}\|_{\vec{V}}^{2}+ 2\| \vec{f}\|_{\vec{V}^{-1}}^{2}.
\end{equation}
 
 Integrating with respect to time, and recalling that $ \overline{\vec{u}_{}}^{m} \in L^{\infty}([0,T],\vec{V}^{\frac{1}{2}}) $ and $\vec{u}_{}^{m}  \in L^{2}([0,T],\vec{V}^{}) $, it follows from Lemma \ref{lemme2} and Lemma	 \ref{l3} that  
  \begin{equation}
\int_{0}^{t}\|\frac{d\vec{u}_{}^{m}}{dt}\|_{\vec{V}^{-1}}^{2}  \le  4( \nu^{2} + C (1/\alpha^{}) K_{2}(T) )^{}
K_{1}(T)+ 2\| \vec{f}\|_{ L^{2}([0,T],\vec{V}^{-1})}^{2}.
\end{equation}

 \subsection{Passing to the limit $\vec{m \rightarrow +\infty}$} From Lemmas \ref{lemme2} and \ref{l4} and  from the reflexivity of the appearing Banach spaces  we can extract a subsequence of $ \vec{u}_{}^{m}$ and $\displaystyle \frac{d\vec{u}_{}^{m}}{dt}$ such that $ \vec{u}_{}^{m}$ converge weakly to some  $\vec{u}_{}$ in $L^{2}([0,T],\vec{V}^{1})$ and $\displaystyle \frac{d\vec{u}_{}^{m}}{dt}$ converge weakly to some $\displaystyle \frac{d\vec{u}_{}^{}}{dt}$ in ${L^{2}([0,T],\vec{V}^{-1})}$
 respectively.\\ 
 Now the interpolation Lemma of Lions and Magenes (Lemma \ref{lionslemma}) implies that 
  $ \vec{u}_{} \in C([0,T];\vec{H})$.\\
 In order to show the convergence of $\vec{u}_{}^{m}$ to $\vec{u}_{}$ in $C([0,T],\vec{H})\cap L^{2}([0,T],\vec{V})$,
we need to show that for $\vec{u}_0 \in  \vec{H} $,  the sequence   $\vec{u}_{}^{m}$ is  a Cauchy sequence in $C([0,T],\vec{H})\cap L^{2}([0,T],\vec{V})$. We know that $\vec{u}_{}^{m}  \in C([0,T],\vec{H})\cap L^{2}([0,T],\vec{V})$.
The difference $\vec{u}_{}^{m+n}  - \vec{u}_{}^{m} $  satifies 
 \begin{equation}
\label{difereence alpha ns}
\begin{array} {llll} \displaystyle
 \frac{d (\vec{u}_{}^{m+n}  - \vec{u}_{}^{m})}{d t}+P_m B( (\overline{\vec{u}_{}}^{m+n} - \overline{\vec{u}_{}}^{m}),\vec{u}_{}^{m}) + P_{m}B(\overline{\vec{u}_{}}^{m+n}, (\vec{u}_{}^{m+n} - \vec{u}_{}^{m}))  \\
 \hskip 4cm + \nu A(\vec{u}_{}^{m+n} - \vec{u}_{}^{m}) =0.
\end{array}
\end{equation}
By taking $\vec{u}_{}^{m+n} - \vec{u}_{}^{m}$ as test function in (\ref{difereence alpha ns}) we get 
 \begin{equation}
\label{d alpha ns}
\begin{array} {llll} \displaystyle
 \displaystyle \frac{1}{2}\frac{d}{d t}\|\vec{u}_{}^{m+n}  - \vec{u}_{}^{m}\|_{\vec{H}}^{2}+ \displaystyle \nu \| \vec{u}_{}^{m+n}- \vec{u}_{}^{m} \|_{\vec{V}}^{2} \\  \hskip 1cm \le
  \displaystyle |(B( (\overline{\vec{u}_{}}^{m+n} - \overline{\vec{u}_{}}^{m}),\vec{u}_{}^{m}),\vec{u}_{}^{m+n} - \vec{u}_{}^{m}) |\\
 \hskip 2cm \displaystyle+ \left\langle P_{m+n}\vec{f}_{\frac{1}{m+n}} - P_m\vec{f}_{\frac{1}{m}}, \vec{u}_{}^{m+n} - \vec{u}_{}^{m} \right\rangle_{\vec{V}^{-1},\vec{V}^{1}}
 \end{array}
 \end{equation}
 Lemma \ref{B:prop 2.3}  combined with Young inequality give
  \begin{equation}
\label{di alpha ns}
\begin{array} {llll} 
\displaystyle |(B( (\overline{\vec{u}_{}}^{m+n} - \overline{\vec{u}_{}}^{m}),\vec{u}_{}^{m}),\vec{u}_{}^{m+n} - \vec{u}_{}^{m}) |
 \le  \displaystyle \frac{1}{\nu}\| \overline{\vec{u}_{}}^{m+n} - \overline{\vec{u}_{}}^{m}\|_{\vec{V}{\frac{1}{2}}}^{2}   \|\vec{u}_{}^{m} \|_{\vec{V}}^{2} \\
 \hskip 4 cm+ \displaystyle\frac{\nu}{4} \|  \vec{u}_{}^{m+n} -  \vec{u}_{}^{m} \|_{\vec{V}}^{2} 
 \end{array}
 \end{equation}
 The duality norm  combined with Young inequality give
  \begin{equation}
\label{di alpha nsforce}
\begin{array} {llll} 
 \displaystyle \left\langle P_{m+n}\vec{f}_{\frac{1}{m+n}} - P_m\vec{f}_{\frac{1}{m}}, \vec{u}_{}^{m+n} - \vec{u}_{}^{m} \right\rangle_{\vec{V}^{-1},\vec{V}^{1}}
 \le  \displaystyle \frac{1}{\nu}\| P_{m+n}\vec{f}_{\frac{1}{m+n}} - P_m\vec{f}_{\frac{1}{m}} \|_{\vec{V}^{-1}}^{2} \\
  \hskip 5cm+ \displaystyle\frac{\nu}{4} \|  \vec{u}_{}^{m+n} -  \vec{u}_{}^{m} \|_{\vec{V}}^{2} 
 \end{array}
 \end{equation}

 Thus we get 
 \begin{equation}
\label{dif alpha ns}
\begin{array} {llll} \displaystyle
 \frac{d}{d t} \|\vec{u}_{}^{m+n} - \vec{u}_{}^{m} \|_{\vec{H}}^{2}+ \nu \| \vec{u}_{}^{m+n} -  \vec{u}_{}^{m} \|_{\vec{V}}^{2}  \le \displaystyle C{\alpha^{-1} \nu^{-1} } \| \vec{u}_{}^{m+n} -\vec{u}_{}^{m}\|_{\vec{H}}^{2}   \|\vec{u}_{}^{m}\|_{\vec{V}}^{2}\\ 
 \hskip 5 cm +\displaystyle \frac{1}{\nu}\| P_{m+n}\vec{f}_{\frac{1}{m+n}} - P_m\vec{f}_{\frac{1}{m}} \|_{\vec{V}^{-1}}^{2}
 \end{array}
 \end{equation}
 By Gr\"{o}nwall inequality we get  
  \begin{equation}
\begin{array} {llll} \displaystyle
 \|\vec{u}_{}^{m+n} - \vec{u}_{}^{m} \|_{\vec{H}}^{2} \\
  \hskip 1cm \le \left(\|\vec{u}_{}^{m+n}(0) - \vec{u}_{}^{m}(0) \|_{\vec{H}}^{2}  +\displaystyle \frac{1}{\nu}\int_{0}^{T}\| P_{m+n}\vec{f}_{\frac{1}{m+n}} - P_m\vec{f}_{\frac{1}{m}} \|_{\vec{V}^{-1}}^{2}  \right)\\
 \hskip 2cm \times \exp{ \displaystyle \int_{0}^{T} \alpha^{-1} \nu^{-1}\|\vec{u}_{}^{m} \|_{\vec{V}}^{2}} dt
   \end{array}
 \end{equation}
 We know that $\vec{u}_{}^{m}  \in  L^{2}([0,T],\vec{V}^{})$, thus there exists $C(\alpha, \nu) \ge 0$ such that $$\exp{\int_{0}^{T} \alpha^{-1} \nu^{-1} \| \vec{u}_{}^{m}\|_{}^{2}dt}  \le C(\alpha,\nu).$$
 We observe that  
 \begin{equation}
\label{dif alpha nsfff}
\begin{array} {llll} \displaystyle
  \int_{0}^{T}\| P_{m+n}\vec{f}_{\frac{1}{m+n}} - P_m\vec{f}_{\frac{1}{m}} \|_{\vec{V}^{-1}}^{2} \le   \displaystyle \underbrace{\int_{0}^{T}\| P_{m+n}\vec{f}_{\frac{1}{m+n}} - P_m\vec{f}_{\frac{1}{m+n}} \|_{\vec{V}^{-1}}^{2}}_{I}\\
  \hskip 5cm+  \displaystyle \underbrace{\int_{0}^{T}\| P_{m}\vec{f}_{\frac{1}{m+n}} - P_m\vec{f}_{\frac{1}{m}} \|_{\vec{V}^{-1}}^{2}}_{II} 
  \end{array}
 \end{equation}
For  $I$ we have
 $$ \| P_{m+n}\vec{f}_{\frac{1}{m+n}} - P_m\vec{f}_{\frac{1}{m+n}} \|_{\vec{V}^{-1}}^{2} \rightarrow 0 $$ because $ \vec{f}_{\frac{1}{m+n}} \in  L^{2}([0,T],\vec{V}^{-1}) $ and 
 $$ \| P_{m+n}\vec{f}_{\frac{1}{m+n}} - P_m\vec{f}_{\frac{1}{m+n}} \|_{\vec{V}^{-1}}^{2} \le \| P_{m+n}\vec{f}_{\frac{1}{m+n}}\|_{\vec{V}^{-1}}^{2} \le \| \vec{f} \|_{\vec{V}^{-1}}^{2} \in L^{1}[0,T]. $$ 
 Using the dominate convergence theorem, we conclude that $I $ tends to zero when $m$ tends to $\infty$.\\
 Similarly by using the dominate convergence theorem combined with the fact that 
 $$   \| P_{m}\vec{f}_{\frac{1}{m+n}} - P_m\vec{f}_{\frac{1}{m}} \|_{\vec{V}^{-1}}^{2}   \rightarrow 0   $$
 and 
 $$   \| P_{m}\vec{f}_{\frac{1}{m+n}} - P_m\vec{f}_{\frac{1}{m}} \|_{\vec{V}^{-1}}^{2} \le \| \vec{f}_{\frac{1}{m+n}} - \vec{f}_{\frac{1}{m}} \|_{\vec{V}^{-1}}^{2} \le \| \vec{f} \|_{\vec{V}^{-1}}^{2} \in L^{1}[0,T],$$
 we obtain that  $II $ tends to zero when $m$ tends to $\infty$.\\
  Since $\vec{u}_0 \in \vec{H}$  then $\|\vec{u}_{}^{m+n}(0) - \vec{u}_{}^{m}(0) \|_{\vec{H}}^{2}$
 converge to zero when $m$ goes to $\infty$.
 By integrating (\ref{dif alpha ns}), we deduce  that $\vec{u}_{}^{m+n} - \vec{u}_{}^{m}$ tend to zero in $C([0,T],\vec{H})\cap L^{2}([0,T],\vec{V})$. 
 This implies that 
    $\vec{u}_{}^{m}$ is  a Cauchy sequence in $C([0,T],\vec{H})\cap L^{2}([0,T],\vec{V})$.\\
    
 Concerning the initial data, we can check that $\vec{u}_{}({0})=\vec{u}_{0}$. Thanks to the result above we have  $ {\vec{u}}_{}^{m}$ converge to $ {\vec{u}}_{}$ in $C([0,T];\vec{H}^{})$, in particular   $\vec{u}_{}^m({0}) $ converge to $ \vec{u}_{}(0)$ in $\vec{H}$. In the other hand, we have that   $\vec{u}_{}^m({0})=P_m \vec{u}_0 $ and  $\vec{u}_{}^m({0})$ converge to $\vec{u}_{0}$ in $\vec{H}$. The unicity of the limit in $\vec{H}$ allows us to deduce the result.\\
 It is obvious from Lemma \ref{regularite} that  $ \overline{\vec{u}}_{}^{m}$ converges to $ \overline{\vec{u}}_{}$ in $C([0,T];\vec{V}^{1/2})$.\\ 
  Moreover, 
 $$ \int_{0}^{T}(A\vec{u}_{}^{m_{}},\phi) dt =-\int_{0}^{T}(A^{1/2}\vec{u}_{}^{m_{}},A^{1/2}\phi)dt \rightarrow -\int_{0}^{T}(A^{1/2}\vec{u}_{}^{},A^{1/2}\phi)dt= \int_{0}^{T}(A\vec{u}_{}^{},\phi) dt,$$
 for all $\phi \in L^{2}([0,T],\vec{V}^{1})$.\\
  Thus, $A\vec{u}_{}^{m_{}}$ converges weakly to  $A\vec{u}_{}^{}$ in $L^{2}([0,T],\vec{V}^{-1})$ as $ m \rightarrow \infty.$\\
 We finish by showing that the non linear term $B(\overline{\vec{u}_{}}^{m_{}} ,\vec{u}_{}^{m_{}} )$ converges weakly to $B(\overline{\vec{u}_{}}^{} ,\vec{u}_{}^{} )$ in $L^{2}([0,T], \vec{V}^{-1})$ as $m \rightarrow \infty.$
 From the properties of the trilinear form we have
 $$ \begin{array}{lllllll}& \displaystyle \left|\int_{0}^{T}(B(\overline{\vec{u}_{}}^{m_{}} ,\vec{u}_{}^{m_{}} ),\phi)dt - \int_{0}^{T}(B(\overline{\vec{u}_{}}^{} ,\vec{u}_{}^{} ),\phi)dt\right|& \\ & \displaystyle \le \int_{0}^{T} \left|(B(\overline{\vec{u}_{}}^{m_{}}-\overline{\vec{u}_{}}^{} ,\vec{u}_{}^{m_{}} ),\phi)\right|dt + \int_{0}^{T} \left|(B(\overline{\vec{u}_{}}^{} ,\vec{u}_{}^{m_{}}- \vec{u}_{}^{}),\phi)\right|dt&
 \end{array},$$
 
 and by using Lemma \ref{B:prop 2.3} combined with H\"{o}lder inequality  we get

$$ \begin{array}{lllllll}  \displaystyle \int_{0}^{T} \left|(B(\overline{\vec{u}_{}}^{m_{}}-\overline{\vec{u}_{}}^{} ,\vec{u}_{}^{m_{}} ),\phi)\right|dt 
  \displaystyle \le  \|\overline{\vec{u}_{}}^{m_{}}-\overline{\vec{u}_{}}^{}\|_{L^{\infty}([0,T],\vec{V}^{1/2})} \|\vec{u}_{}^{m_{}}\|_{L^{2}([0,T],\vec{V}^{1})} \|\phi\|_{L^{2}([0,T],\vec{V}^{1})},
\end{array}$$
and
$$ \begin{array}{lllllll}  \displaystyle \int_{0}^{T} \left|(B(\overline{\vec{u}_{}}^{} ,\vec{u}_{}^{m_{}}- \vec{u}_{}^{}),\phi)\right|dt
  \displaystyle \le   \displaystyle  \|\overline{\vec{u}_{}}^{} \|_{L^{\infty}([0,T],\vec{V}^{1/2})}\|\vec{u}_{}^{m_{}}- \vec{u}_{}^{}\|_{L^{2}([0,T],\vec{V}^{1})}\|\phi\|_{L^{2}([0,T],\vec{V}^{1})}.\\
\end{array}$$
 

Thus $B(\overline{\vec{u}_{}}^{m_{}} ,\vec{u}_{}^{m_{}} )$ converges weakly to $B(\overline{\vec{u}_{}}^{} ,\vec{u}_{}^{} )$ in $L^{2}(([0,T],\vec{V}^{-1}))$ as $m \rightarrow \infty.$ This implies  that $P_m B(\overline{\vec{u}_{}}^{m_{}} ,\vec{u}_{}^{m_{}} )$ converges weakly to $B(\overline{\vec{u}_{}}^{} ,\vec{u}_{}^{} )$ in $L^{2}(([0,T],\vec{V}^{-1}))$ as $m \rightarrow \infty.$\\
The convergence of $P_m \vec{f}_{\frac{1}{m}}$ to $\vec{f}$ in $L^{2}([0,T],\vec{V}^{-1})$ is obvious because $\vec{f} \in L^{2}([0,T],\vec{V}^{-1}).$    We have shown that  $\vec{u}_{}$ satisfies (\ref{alpha ns prime}) viewed as a functional equality in $\vec{V}^{-1}$.\\

\subsection{Uniqueness}

The solution constructed above is unique. Next, we will show the continuous dependence of the  solutions on the initial data and in particular the uniqueness.\\
Let $\vec{u}_{}$ and $\vec{v}_{}$ be any two solutions of (\ref{alpha ns prime}) on the interval $[0,T]$, with initial values  $\vec{u}_{0} \in \vec{H}^{} $ and $\vec{v}_{0} \in \vec{H}^{} $, respectively. Let us denote by  $\vec{w}_{} =\vec{u}_{}-\vec{v}_{}$ and $\overline{\vec{w}_{}} =\overline{\vec{u}_{}}-\overline{\vec{v}_{}}$.
 Then,  we can write the evolution equation for $\vec{w}_{}$ as an equality in $\vec{V}^{-1}$ given by 

 \begin{equation}
 \label{evolution1}
\frac{d }{d t}\vec{w}_{}+ \nu A \vec{w}_{} +B(\overline{\vec{w}_{}},\vec{u}_{}) +B(\overline{\vec{v}},\vec{w}_{}) = 0
\end{equation}

We take the inner product of (\ref{evolution1})  with $\vec{w}^{}$.  Applying the Lemma \ref{lionslemma} of Lions-Magenes and by the properties of the trilinear form  we get 
\begin{equation}
 \frac{1}{2} \frac{d}{dt}
\|{\vec{w}}_{}\|_{\vec{H}}^{2} + \nu \|\nabla \vec{w}_{}\|_{\vec{H}}^{2} +( B(\overline{\vec{w}}_{},\vec{u}_{}),\vec{w}_{})=0
\end{equation}
Using  Lemma \ref{B:prop 2.3} combined with Young's inequality and Lemma \ref{regularite}, we obtain

 
  \begin{equation}
 \label{jfkjfkjkdv}\begin{array}{llllll}
 \displaystyle \frac{d}{dt}
\|{\vec{w}}_{}\|_{\vec{H}}^{2} +\nu \|\nabla
{\vec{w}}_{}\|_{\vec{H}}^{2}  &\le& \displaystyle \frac{1}{\alpha}\frac{4}{\nu} \|{\vec{w}}\|_{\vec{H}}^{2}
\|{\vec{u}}_{} \|^{2}_{\vec{V}} 
\end{array}
\end{equation}

Using Gr\"{o}nwall's inequality we conclude that
$$ \|{\vec{w}}_{}\|_{L^{\infty}([0,T],\vec{H}^{})}^{2} \le \|{\vec{w}}^{}_{0}\|_{\vec{H}}^{2} \exp{\frac{1}{\alpha}{\frac{4}{\nu}
K_{1}}}.$$ 
We have shown the continuous dependence of the solutions on the inital data in the $L^{\infty}([0,T],\vec{H})$  norm. In particular, if ${\vec{w}}^{}_{0}=0$ then ${\vec{w}}_{}=0$ and the solutions are unique for all $t \in [0,T] .$\\ 
\begin{Remark}
Since $T>0$ is arbitrary, the  solution above may be uniquely extended for all time.
\end{Remark}
\begin{Remark}
\label{prague}
The pressure is absent in equations (\ref{alpha ns prime}) once we can recover it by using the de Rham Theorem \cite{LuTa78}.
Generally, the existence of the pressure  is not obvious and the pressure may not exist \cite{27}.\\
In our situation, as we work with   periodic boundary conditions, we can proceed differently. 
 We take the divergence of
(\ref{alpha ns}), this yields   the following equation for the pressure
\begin{equation}
\label{elliptic pressure} -\Delta p_{} = \nabla \cdot (\nabla \cdot
\overline{\vec{u}_{}}\vec{u}_{}),
\end{equation}
where $\overline{\vec{u}}\vec{u}$ is the tensor $(\overline{u}^{i} u^{j})_{1 \le i,j \le 3}$.\\
We recall that for any divergence-free field $\vec{u}$, we have $$
(\overline{\vec{u}}\cdot\nabla)\vec{u}=\nabla\cdot(\overline{\vec{u}}\vec{u}).$$
Since $(\overline{\vec{u}_{}}\cdot\nabla)\vec{u}_{}   \in L^{2}([0,T],\vec{H}^{-1})$.
We  conclude from  the classical elliptic theory that $p_{}$ is
bounded in $\displaystyle L^{2}([0,T],L^2(\tore)).$ 
\end{Remark}
\begin{lemme}
\label{Uniqueness of the pressure} (Uniqueness of the pressure) The pressure is uniquely determined by
the velocity.
\end{lemme}
\textbf{Proof.}
Let $({\vec{u}_{},p_{}})$ and $({\vec{u}_{},q_{}})$ be two weak solutions of (1.1) (with $\theta \ge \frac{1}{4}$) such that $\int_{\tore}p_{} = \int_{\tore}q_{}=0$.  Then  for every $ {\phi} \in \vec{H}^1$ and for almost every $ t \in \RR^{+}$ we have 
$$(\nabla(p_{}-q_{}), {\phi})=0,$$
 combined with the condition $\displaystyle \int_{\tore}p_{}=  \displaystyle \int_{\tore}q_{} =0$  gives $\displaystyle p_{} = q_{}$.\\
This concludes the proof of Theorem \ref{2}. 

\section{Proof of Theorem \ref{1}}

\subsection{ Energy estimates in $\vec{V}^{\frac{1}{2}}$ }
Let us write the energy inequality in the space $\vec{V}^{\frac{1}{2}}$.  Taking the $\vec{V}^{\frac{1}{2}}$ scalar product of eqs. (\ref{galerkine}) with  $\vec{u}_{}^{m}$ it turns out, due to the divergence-free condition, that 
\begin{equation}
 \frac{1}{2} \frac{d}{dt}
\|{\vec{u}}_{}^{m}\|_{\frac{1}{2}}^{2} + \nu \|\nabla \vec{u}_{}^{m}\|_{\vec{V}^{\frac{1}{2}}}^{2} =( 
P_{m} B(\overline{\vec{u}_{}}^{m},\vec{u}^{m}_{}), \vec{u}_{}^{m})_{\vec{V}^\frac{1}{2}} + (P_{m}f_{\frac{1}{m}}, \vec{u}_{}^{m})_{\vec{V}^\frac{1}{2}}.
\end{equation}
By definition of the scalar product on $\vec{V}^{\frac{1}{2}}$, we get 
\begin{equation}
\left|(P_{m} B(\overline{\vec{u}_{}}^{m},\vec{u}^{m}_{}), \vec{u}_{}^{m})_{\vec{V}^\frac{1}{2}}\right| \le  \| B(\overline{\vec{u}_{}}^{m},\vec{u}^{m}_{})\|_{\vec{V}^{-{\frac{1}{2}}}} \|\nabla \vec{u}_{}^{m} \|_{\vec{V}^{\frac{1}{2}}}.
\end{equation}

The Sobolev embedding in Theorem \ref{soblevin} together with a H\"{o}lder estimate give
\begin{equation}
\begin{array}{llllll}
  \left|(  B(\overline{\vec{u}_{}}^{m},\vec{u}^{m}_{}) ,|\nabla| \vec{u}_{}^{m})_{L^2}\right|&\le & c\|B(\overline{\vec{u}_{}}^{m},\vec{u}^{m}_{}) \|_{\vec{V}^{-1/2}} \|\nabla \vec{u}_{}^{m} \|_{\vec{V}^{1/2}} \\ 
  &\le& c\| B(\overline{\vec{u}_{}}^{m},\vec{u}^{m}_{})\|_{{L}^{3/2}} \|\nabla \vec{u}_{}^{m} \|_{\vec{V}^{1/2}} \\ 
   &\le& c \| \overline{\vec{u}_{}}^{m}\|_{{L}^{6}} \|\nabla\vec{u}_{}^{m}\|_{L^2} \|\nabla \vec{u}_{}^{m} \|_{\vec{V}^{1/2}} \\
    &\le& c \| \overline{\vec{u}_{}}^{m}\|_{\vec{V}}\|\vec{u}_{}^{m}\|_{\vec{V}} \|\nabla \vec{u}_{}^{m} \|_{\vec{V}^{1/2}}. 
  \end{array}
\end{equation}

By the interpolation inequality between $\vec{V}^{1/2}$ and $\vec{V}^{3/2}$ we get 
\begin{equation}
\begin{array}{llllll}
  \left|(B(\overline{\vec{u}_{}}^{m},\vec{u}^{m}_{}) ,|\nabla| \vec{u}_{}^{m})_{L^2} \right|&\le& \| \overline{\vec{u}_{}}^{m}\|_{\vec{V}}\|\vec{u}_{}^{m}\|_{\vec{V}^{1/2}}^{1/2} \|\nabla \vec{u}_{}^{m} \|_{\vec{V}^{1/2}}^{3/2}. 
  \end{array}
\end{equation}

Using the fact that 
\begin{equation}
\begin{array}{llllll}
 (P_{m}\vec{f}_{\frac{1}{m}},  \vec{u}_{}^{m})_{\vec{V}^{1/2}} \le \| \vec{f}  \|_{\vec{V}^{-1/2}}\|\nabla \vec{u}_{}^{m} \|_{\vec{V}^{1/2}}
 \end{array}
\end{equation}
combined with the Young inequality, we infer 

\begin{equation}
\label{energie}
 \frac{d}{dt}
\|{\vec{u}}_{}^{m}\|_{\vec{V}^{1/2}}^{2} + \nu \|\nabla \vec{u}_{}^{m}\|_{\vec{V}^{1/2}}^{2} \le \frac{9}{2 \nu^{3}}\| \overline{\vec{u}}^{m}\|_{\vec{V}}^{4}\|\vec{u}_{}^{m}\|_{\vec{V}{1/2}}^{2} +\frac{2}{\nu}\| \vec{f}  \|_{\vec{V}^{-1/2}}^{2}.
\end{equation}

The following Lemmas are crucial to the proof of Theorem 1.1.

\begin{lemme}
\label{3.3}
Let $\vec{f} \in L^{2}([0,T],\vec{V}^{-{\frac{1}{2}}})$ and $\vec{u}_{0} \in \vec{V}^{\frac{1}{2}}$, there exists $M_{1}(T)$ and $M_{2}(T)$ independent of m such that the solution $\vec{u}_{}^{m}$ to the Galerkin truncation (\ref{galerkine}) satisfies 
\begin{equation}
\|\vec{u}_{}^{m}\|_{L^{2}([0,T],\vec{V}^{3/2})}^{2} \le M_{1}(T), \ \hbox{for all} \  T \ge 0,
\end{equation}
and
\begin{equation}
\|\vec{u}_{}^{m}\|_{L^{\infty}([0,T],\vec{V}^{1/2})}^{2} \le M_{2}(T), \  \hbox{for all} \ t \ge 0.
\end{equation}
\end{lemme}   


\textbf{Proof.} 
We proceed by two steps:\\
\underline{First step:}
From (\ref{energie}) and by using  Gr\"{o}nwall's inequality we conclude that 
\begin{equation}
\|\vec{u}_{}^{m}\|_{L^{\infty}([0,T],\vec{V}^{1/2})}^{2} \le M_{2}(T), 
\end{equation}
{where} $M_{2}(T)$ is given by 
$$M_{2}(T)= \left(\|{\vec{u}}^{}_{0}\|_{\vec{V}^{1/2}}^{2} 
 +\frac{2}{\nu}\|\vec{f}\|_{L^{2}([0,T],\vec{V}^{-1/2})}^{2}\right)\exp{\frac{9}{2\nu^{3}\alpha^{}} K_{1}(T)^{1/2} K_{2}(T)^{1/2}}.$$ 

\underline{Second step:}
Integrate the original inequality  (\ref{energie}) on $[0,T]$.
We get 
\begin{equation}
\begin{array}{lllll}
\|{\vec{u}}_{}^m(T,\vec{x})\|_{\vec{V}^{1/2}}^{2} +\nu \displaystyle \int_{0}^{T}\|\nabla
{\vec{u}}_{}^m(t,\vec{x})\|_{\vec{V}^{1/2}}^{2} dt &\le \displaystyle \|{\vec{u}}^{}_{0}\|_{\vec{V}^{1/2}}^{2}
 +\frac{2}{\nu}\| \vec{f}\|_{L^{2}([0,T],\vec{V} ^{-1/2})}^{2}& \\
& \displaystyle + \frac{9}{2\nu^{3}\alpha^{}} M_{2}(T) K_{1}(T)^{1/2} K_{2}(T)^{1/2}.&
\end{array}
\end{equation}
We set $$M_{1}(T)=\displaystyle \|{\vec{u}}^{}_{0}\|_{\vec{V}^{1/2}}^{2}
 +\frac{2}{\nu}\| \vec{f}\|_{L^{2}([0,T],\vec{V} ^{-1/2})}^{2}
\displaystyle + \frac{9}{2\nu^{3}\alpha^{}} M_{2}(T) K_{1}(T)^{1/2} K_{2}(T)^{1/2}.$$

Thus $\vec{u}_{}^{m} \in  L^{\infty }([0,T],\vec{V}^{1/2}) \cap  L^{2}([0,T],\vec{V}^{3/2})$ for all $T>0$.\\

We deduce from Lemma \ref{regularite} that 
\begin{equation}
\|\overline{\vec{u}_{}}^{m}\|_{L^{2}([0,T],\vec{V}^{2})}^{2} \le( 1/ \alpha^{}) M_{1}(T), \ \hbox{for all} \ T \ge 0 
\end{equation}
and
\begin{equation}
\|\overline{\vec{u}_{}}^{m}\|_{L^{\infty}([0,T],\vec{V}^{1})}^{2} \le (1/\alpha^{}) M_{2}(T), \ \hbox{for all} \ T \ge 0   
\end{equation}

Our next result provides an estimate on the time derivative of $\vec{u}_{}^{m}$.
\begin{lemme}
\label{3.5}
Let $\vec{f} \in L^{2}([0,T],\vec{V}^{-1/2})$ and $\vec{u}_{0} \in \vec{V}^{1/2}$, there exists $M_{3}(T)$ independent of $m$ such that the time derivative of the solution $\vec{u}_{}^{m}$ to the Galerkin truncation \ref{galerkine} satisfies 
\begin{equation}
\label{du sur dt}
\begin{array}{llll}
\displaystyle
\|\frac{d\vec{u}_{}^{m}}{dt}\|_{L^{2}([0,T],\vec{V}^{-1/2})}^2 \le M_{3}(T), \ \hbox{for all} \ T \ge 0   , \\
\end{array}
\end{equation}
where 
\begin{equation}
\displaystyle M_{3}(T)= \left(  4( \nu^{2} + C  M_{1}(T) )^{}
M_{2}(T)+ 2\| \vec{f}\|_{L^{2}([0,T],\vec{V} ^{-1/2}) }^{2}    \right). 
\end{equation}

\end{lemme}
\textbf{Proof.}
Taking  the $\vec{V}^{-1/2}$ norm of  (\ref{galerkine}), one obtains:
\begin{equation}
\|\frac{d\vec{u}_{}^{m}}{dt}\|_{\vec{V}^{-1/2}}  \le  \nu
\|{\vec{u}_{}^{m}}{}\|_{\vec{V}^{1/2}}+\| B(\overline{\vec{u}_{}}^{m} ,\vec{u}_{}^{m} )\|_{\vec{V}^{-1/2}}+ \| \vec{f}\|_{\vec{V}^{-1/2} }
\end{equation}
where we note that 
$$\| B(\overline{\vec{u}_{}}^{m} ,\vec{u}_{}^{m} )\|_{\vec{V}^{-1/2}}= \sup{\left\{(B(\overline{\vec{u}_{}}^{m} ,\vec{u}_{}^{m} ),w); w \in \vec{V}^{1/2},  \|w\|_{{1/2}} \le 1 \right\}}.$$

It remain to show that $\| B(\overline{\vec{u}_{}}^{m} ,\vec{u}_{}^{m} )\|_{\vec{V}^{-1/2}}$  is bounded.\\
Indeed, we have from Lemma \ref{B:prop 2.3} that


\begin{equation}
\begin{array}{cccccccc}
  \|B(\overline{\vec{u}_{}}^{m} ,\vec{u}_{}^{m} )\|_{\vec{V}^{-1/2} }
 &\le& C \|\overline{\vec{u}_{}}^{m} \|_{\vec{V}^{1/2}} \|\vec{u}_{}^{m}\|_{\vec{V}^{3/2}}\\
 &\le& C \|{\vec{u}}^{m}_{} \|_{\vec{V}^{1/2}} \|\vec{u}_{}^{m}\|_{\vec{V}^{3/2}}.
 \end{array}
 \end{equation}

 It follows that 
 \begin{equation}
\|\frac{d\vec{u}_{}^{m}}{dt}\|_{\vec{V}^{-1/2}}^{2}  \le  4( \nu^{2} + C \|{\vec{u}}^{m} \|_{\vec{V}^{3/2}}^{2})^{}
\|{\vec{u}_{}^{m}}{}\|_{\vec{V}^{1/2}}^{2}+ 2\| \vec{f}\|_{\vec{V}^{-1/2} }^{2}.
\end{equation}
 
 Integrating with respect to time, it follows from Lemma \ref{3.3}  that  
  \begin{equation}
\int_{0}^{T}\|\frac{d\vec{u}_{}^{m}}{dt}\|_{\vec{V}^{-1/2}}^{2}  \le  4( \nu^{2} + C  M_{1}(T) )^{}
M_{2}(T)+ 2\| \vec{f}\|_{L^{2}([0,T], \vec{V}^{-1/2}) }^{2}
\end{equation}
yields (\ref{du sur dt}) and the proof is complete.\\

\subsection{Passing to the limit}
 We are now ready to take limits along subsequences of $ \vec{u}_{}^{m}$ to show the existence of solutions that satisfy (\ref{alpha ns prime}) viewed as a functional equality in $\vec{V}^{-\frac{1}{2}}$.\\
  From Lemmas \ref{3.3} and \ref{3.5} combined with Alaoglu compactness theorem we can extract a subsequence of $ \vec{u}_{}^{m}$ and $\displaystyle \frac{d\vec{u}_{}^{m}}{dt}$ that converge weakly to some  $\vec{u}_{}$ in $L^{2}([0,T],\vec{V}^{\frac{3}{2}})$ and $\displaystyle \frac{d\vec{u}_{}^{}}{dt}$ in ${L^{2}([0,T],\vec{V}^{-\frac{1}{2}})}$
 respectively.\\ 
 Now the  Lions-Magenes Lemma \ref{lionslemma}, implies that $ \vec{u}_{} \in C([0,T];\vec{V}^{\frac{1}{2}})$ and this subsequence  $ \vec{u}_{}^{m}$ converge to $ \vec{u}_{}$ in $C([0,T];\vec{H})$ and in particular we have $\vec{u}_{}({0})=\vec{u}_{0}$.
 Thus is obvious that  $ \overline{\vec{u}_{}}^{m}$ converge to $ \overline{\vec{u}_{}}$ in $C([0,T];\vec{V}^{\frac{1}{2}})$.\\ 
 Fix $\epsilon$ such that $0< \epsilon < 1/2$, Since $\vec{V}^{\frac{3}{2}}$ is compactly embedded in $\vec{V}^{\frac{3}{2}-\epsilon}$ then Theorem \ref{compact},  
  implies that there exists a further subsequence of $ \vec{u}_{}^{m}$,  denoted by $ \vec{u}_{}^{m_{i}}$, that converges strongly in $L^{2}([0,T];\vec{V}^{\frac{3}{2}-\epsilon})$. Moreover, 
 $$ \int_{0}^{T}(A\vec{u}_{}^{m_{i}},\phi) dt =\int_{0}^{T}(\vec{u}_{}^{m_{i}},A^{}\phi)dt \rightarrow \int_{0}^{T}(\vec{u}_{}^{},A^{}\phi)dt= \int_{0}^{T}(A\vec{u}_{}^{},\phi) dt,$$
 for all $\phi \in L^{2}([0,T],\vec{V}^{\frac{1}{2}})$. Thus, $A\vec{u}_{}^{m_{i}}$ converge weakly to  $A\vec{u}_{}^{}$ in $L^{2}([0,T],\vec{V}^{-\frac{1}{2}})$ as $ i \rightarrow \infty.$\\
 We finish by showing that the non linear term $B(\overline{\vec{u}_{}}^{m_{i}} ,\vec{u}_{}^{m_{i}} )$ converge weakly to $B(\overline{\vec{u}_{}}^{} ,\vec{u}_{}^{} )$ in $L^{2}([0,T],\vec{V}^{-\frac{1}{2}})$ as $i \rightarrow \infty.$
 From the properties of the trilinear form we have
 $$ \begin{array}{lllllll}& \displaystyle \left|\int_{0}^{T}(B(\overline{\vec{u}_{}}^{m_{i}} ,\vec{u}_{}^{m_{i}} ),\phi)dt - \int_{0}^{T}(B(\overline{\vec{u}_{}}^{} ,\vec{u}_{}^{} ),\phi)dt\right|& \\ & \displaystyle \le \int_{0}^{T} \left|(B(\overline{\vec{u}_{}}^{m_{i}}-\overline{\vec{u}_{}}^{} ,\vec{u}_{}^{m_{i}} ),\phi)\right|dt + \int_{0}^{T} \left|(B(\overline{\vec{u}_{}}^{} ,\vec{u}_{}^{m_{i}}- \vec{u}_{}^{}),\phi)\right|dt&
 \end{array},$$
 by using H\"{o}lder inequality combined with Sobolev injection Theorem and Poincare inequality  we get 
$$ \begin{array}{llllllllll} & \displaystyle \int_{0}^{T} \left|(B(\overline{\vec{u}_{}}^{m_{i}}-\overline{\vec{u}_{}}^{} ,\vec{u}_{}^{m_{i}} ),\phi)\right|dt + \int_{0}^{T} \left|(B(\overline{\vec{u}_{}}^{} ,\vec{u}_{}^{m_{i}}- \vec{u}_{}^{}),\phi)\right|dt&\\
  &\displaystyle \le \int_{0}^{T} \|\overline{\vec{u}_{}}^{m_{i}}-\overline{\vec{u}_{}}^{}\|_{\vec{V}^{\frac{1}{2}}} \|\vec{u}_{}^{m_{i}}\|_{\vec{V}^{\frac{3}{2}}} \|\phi\|_{\vec{V}^{\frac{1}{2}}}dt + \int_{0}^{T} \|\overline{\vec{u}_{}}^{} \|_{\vec{V}^{\frac{1}{2}+\epsilon}}\|\vec{u}_{}^{m_{i}}- \vec{u}_{}^{}\|_{\vec{V}^{\frac{3}{2}-\epsilon}}\|\phi\|_{\vec{V}^{\frac{1}{2}}} dt&\\
  &\displaystyle \le \displaystyle\|\phi\|_{L^{2}([0,T],\vec{V}^{\frac{1}{2}})} \|\overline{\vec{u}_{}}^{m_{i}}-\overline{\vec{u}_{}}^{}\|_{L^{\infty}([0,T],\vec{V}^{\frac{1}{2}})}\|\vec{u}_{}^{m_{i}}\|_{L^{2}([0,T],\vec{V}^{\frac{3}{2}})}& \\
  & \hskip 2cm \displaystyle + \displaystyle  c_{p} \|\phi\|_{L^{2}([0,T],\vec{V}^{\frac{1}{2}})} \|\overline{\vec{u}_{}}^{} \|_{L^{\infty}([0,T],\vec{V}^{1})}\|\vec{u}_{}^{m_{i}}- \vec{u}_{}^{}\|_{L^{2}([0,T],\vec{V}^{\frac{3}{2}-\epsilon})}. &
 \end{array}$$
 Thus $B(\overline{\vec{u}_{}}^{m_{i}} ,\vec{u}_{}^{m_{i}} )$ converge weakly to $B(\overline{\vec{u}_{}}^{} ,\vec{u}_{}^{} )$ in $L^{2}([0,T],\vec{V}^{-1/2})$ as $i \rightarrow \infty.$\\
  We have shown that  $\vec{u}_{}$ satisfies (\ref{alpha ns prime}) viewed as a functionnal equality in $\vec{V}^{-1/2}$.\\

\subsection{Uniqueness}
 The solution constructed above is unique. Next, we will show the continuous dependence of the  solutions on the initial data and in particular the uniqueness.\\
 Let $\vec{u}_{}$ and $\vec{v}_{}$ any two solutions of (\ref{alpha ns}) on the interval $[0,T]$, with initial values  $\vec{u}_{0} \in \vec{V}^{\frac{1}{2}} $ and $\vec{v}_{0} \in \vec{V}^{\frac{1}{2}} $, respectively. Let us denote by  $\vec{w}_{} =\vec{u}_{}-\vec{v}_{}$ and $\overline{\vec{w}_{}} =\overline{\vec{u}_{}}-\overline{\vec{v}_{}}$.
 Then from (\ref{alpha ns prime}) we can write the evolution equation for $\vec{w}_{}$ as an equality in $\vec{V}^{-\frac{1}{2}}$ given by 
\begin{equation}
 \label{difference}
\frac{d }{d t}\vec{w}_{}+ \nu A \vec{w}_{} +B(\overline{\vec{w}_{}},\vec{u}_{}) +B(\overline{\vec{v}_{}},\vec{w}_{}) = 0
\end{equation}

We take the inner product of \ref{difference}  with $\vec{w}^{}$,  applying the Lemma of Lions-Magenes \ref{lionslemma} and by the 
 bilinearity of $B$  we have
\begin{equation}
 \frac{1}{2} \frac{d}{dt}
\|{\vec{w}}_{}\|_{\vec{V}^{1/2}}^{2} + \nu \|\nabla \vec{w}_{}\|_{\vec{V}^{1/2}}^{2} +( B(\overline{\vec{w}_{}},\vec{u}_{}),| \nabla|\vec{w}_{} )_{L^2}+( B(\overline{\vec{v}_{}},\vec{w}_{}),| \nabla|\vec{w}_{} )_{L^2} =0
\end{equation}

The first non linear term is estimated by 
\begin{equation}
  \begin{array}{cccc}
 |( B(\overline{\vec{w}_{}},\vec{u}_{}),|\nabla|\vec{w}_{})_{L^2}| &\le&  \|\overline{\vec{w}}
\nabla \vec{u}_{} \|_{{\vec{V}^{-1/2}}} \| \nabla \vec{w}_{}\|_{\vec{V}^{1/2}}\\
&\le&  C\|\overline{\vec{w}_{}}\|_{\vec{V}^{1/2}}
\|\nabla \vec{u}_{} \|_{\vec{V}^{1/2}} \| \nabla \vec{w}_{}\|_{\vec{V}^{1/2}},\\ 
\end{array}
\end{equation}
and by Young's inequality, we obtain
 \begin{equation}
|( B(\overline{\vec{w}_{}},\vec{u}_{}),|\nabla|\vec{w}_{})_{L^2}| \le \frac{C}{\nu}  \|\overline{\vec{w}}\|_{\vec{V}^{1/2}}^2
\| \vec{u}_{} \|_{\vec{V}^{3/2}}^{2}+ \frac{\nu}{4}\| \nabla \vec{w}^{}\|_{\vec{V}^{1/2}}^{2}.
\end{equation}

The second non linear term is estimated by 
\begin{equation} \begin{array}{cccc}|(B(\overline{\vec{v}_{}},\vec{w}_{}),|\nabla|\vec{w}_{})_{L^2}| &\le& \|\overline{\vec{v}_{}}
\nabla \vec{w}_{} \|_{\vec{H}} \| \nabla \vec{w}_{}\|_{\vec{H}}\\
&\le&  C\|\overline{\vec{v}_{}}\|_{\vec{V}} \|\nabla \vec{w}_{} \|_{\vec{V}^{1/2}} \| \nabla \vec{w}_{}\|_{\vec{H}}\\
 &\le&  C\|{\vec{v}_{}}\|_{\vec{V}} \|\vec{w}_{} \|_{\vec{V}^{3/2}}^{3/2} \|  \vec{w}_{}\|_{\vec{V}^{1/2}}^{1/2}\\
\end{array}
\end{equation}
and by Young's inequality, we obtain
 \begin{equation}
|( B(\overline{\vec{v}_{}}_{},\vec{w}_{}),|\nabla| \vec{w}_{})_{L^2}| \le \frac{C}{\nu}   \|{\vec{v}}_{}\|_{\vec{V}}^{4} \|  \vec{w}_{}\|_{\vec{V}^{1/2}}^{2}  + \frac{\nu}{4}\| \nabla \vec{w}^{}\|_{\vec{V}^{1/2}}^{2}.
\end{equation}

From the above inequalities we get
 \begin{equation}
 \label{jfkjfkjkdvprime}\begin{array}{llllll}
  \displaystyle \frac{d}{dt}
\|{\vec{w}}_{}\|_{\vec{V}^{1/2}}^{2} +\nu \|\nabla
{\vec{w}}_{}\|_{\vec{V}^{1/2}}^{2}  &\le& \displaystyle
\frac{C}{\nu}  \|  \vec{w}_{}\|_{\vec{V}^{1/2}}^{2} \left( \|{\vec{v}_{}}\|_{\vec{V}}^{4}+ \|\vec{u}_{} \|_{\vec{V}^{3/2}}^{2}  \right)
\end{array}
\end{equation}
Using Gr\"{o}nwall's inequality we conclude that
$$ \|{\vec{w}}_{}\|_{L^{\infty}([0,T],\vec{V}^{1/2})}^{2} \le \|{\vec{w}}^{}_{0}\|_{\vec{V}^{1/2}}^{2} \exp{{\frac{C}{\nu} \left( \int^{T}_{0}\|{\vec{v}_{}}\|_{\vec{V}}^{4}dt+\int^{T}_{0} \|\vec{u}_{} \|_{\vec{V}^{3/2}}^{2}dt  \right)}}.$$ 
Note  that the above estimate does not depend on $\alpha$, thus we have not used the regularization effect to get the uniqueness in $\vec{V}^{\frac{1}{2}}$ and this is not surprising because  the solutions of the Navier-Stokes equations with initial data in $\vec{V}^{\frac{1}{2}}$ are unique.\\
We have shown the continuous dependence of the solutions on the inital data in the $L^{\infty}([0,T],\vec{V}^{1/2})$  norm. In particular, if ${\vec{w}}^{}_{0}=0$ then ${\vec{w}}_{}=0$ and solutions are unique for all $t \in [0,T] .$ Since $T>0$ is arbitrary this solution may be uniquely extended for all time.\\ 

\begin{Remark}
The pressure is absent in eqs. (\ref{alpha ns prime}) once we can reconstruct it from $\overline{\vec{u}_{}}$ and $\vec{u}_{} $ (up to a constant), if necessary, see Lemma \ref{Uniqueness of the pressure} and Remark \ref{prague}.
 We take the divergence of
(\ref{alpha ns}). This yields the following equation for the pressure
\begin{equation}
 -\Delta p_{} = \nabla \cdot (\nabla \cdot
\overline{\vec{u}_{}}\vec{u}_{}).
\end{equation}
Where $\overline{\vec{u}}\vec{u}$ is the tensor $(\overline{u}^{i} u^{j})_{1 \le
i,j \le 3}$.\\
We recall that for any field $\vec{u}$ with divergence is equal
to zero, we have $$
(\overline{\vec{u}}\cdot\nabla)\vec{u}=\nabla\cdot(\overline{\vec{u}}\vec{u}).$$
Since $(\overline{\vec{u}_{}}\cdot\nabla)\vec{u}_{}   \in L^{2}([0,T],\vec{H}^{-1/2})$.
One concludes from  the classical elliptic theory that $p_{}$ is
bounded in $\displaystyle L^{2}([0,T],{H}^{1/2}(\tore)).$ 
\end{Remark}

This finish the proof of Theorem \ref{1}.

\section{Relation between NSE and Leray-$\alpha$: construction of suitable solutions}
 The regularized solution contructed above with $\theta \ge \frac{1}{4}$ is unique. So the solution  is suitable. 
 In this section we will construct a suitable weak solution to the Navier-Stokes equations by taking the limit when $\alpha$ tends to zero to the regularized solution.

We start by the following Lemma: 
\begin{lemme}
\label{lemmelocalen}
Let $(\vec{u}_{\alpha}, p_{\alpha})$ be the unique solution of  (\ref{alpha ns}) with  $\theta \ge \frac{1}{4}$ then  $(\vec{u}_{\alpha}, p_{\alpha})$ verifies the following local inequality
\begin{equation}
\label{local alpha}
\begin{array}{lcc}
\displaystyle 2\nu\int_{0}^{T}\int_{\tore}|\nabla \vec{u}_{\alpha}|^{2}\phi \ d\vec{x}dt\\
\hskip 2cm= \displaystyle \int_{0}^{T}\int_{\tore}| \vec{u}_{\alpha}|^{2}\left(\phi_t + \nu \Delta \phi \right) + \left(| \vec{u}_{\alpha}|^{2} \overline{\vec{u}_{\alpha}}+ 2p\vec{u}_{\alpha}\right) \cdot \nabla \phi   \  d\vec{x}dt\\
 \hskip 4cm \displaystyle + \displaystyle 2\int_{0}^{T}\int_{\tore} \vec{f} \vec{u}\phi \  d\vec{x}dt, 
\end{array}
\end{equation}
for all $T \in (0,+\infty]$ and for all  non negative fonction  $\phi \in C^{\infty}$ and  supp $\phi \subset\subset \tore\times(0,T).$
\end{lemme}

\textbf{Proof.} 
We take $2\phi\vec{u}_{\alpha}$ as test function in (\ref{alpha ns}). We note that the condition $\theta \ge 1/4 $  ensure that all the tems are well defined. In particular the integral $$\displaystyle 2\int_{0}^{T} \int_{\tore} \overline{\vec{u}_{\alpha}}\nabla \vec{u}_{\alpha} \cdot \vec{u}_{\alpha} \phi   \  d\vec{x} dt $$ is finite by using the fact that  $ \overline{\vec{u}_{\alpha}}\nabla \vec{u}_{\alpha} \in L^2([0,T];\vec{H}^{-1})$ and $2\phi\vec{u}_{\alpha} \in L^2([0,T];\vec{H}^{1})$.\\
An integration by part combined with   $\phi(T,\cdot)=\phi(0,\cdot)=0$  and the following identity 
\begin{equation}
2\int_{\tore} \overline{\vec{u}_{\alpha}}\nabla \vec{u}_{\alpha} \cdot \vec{u}_{\alpha} \phi   \  d\vec{x} = \int_{\tore} \overline{\vec{u}_{\alpha}} |\vec{u}_{\alpha}|^2   \cdot \nabla \phi   \  d\vec{x}
\end{equation}
yield that  $(\vec{u}_{\alpha}, p_{\alpha})$ verifies  (\ref{local alpha}).\\

In order to take the limits  $\alpha \rightarrow 0$ over  (\ref{local alpha}), we need first to show that for all $\bu_\alpha \in L^p(0,T;L^p(\tore)^3) $, $2 \le p<10/3,$  we have:  \begin{align}
\overline{\bu_\alpha}  &\rightarrow \bu &&\textrm{strongly in  }
L^p(0,T;L^p(\tore)^3) \textrm{ for all }  2 \le p<10/3.
\end{align}
This is the aim of the two following Lemmas. 
%
\begin{lemme}
\label{fourierdiscret}
Let $ \theta \in\RR^+,  0 \le \beta \le 2\theta , s\in\RR$ and assume that $\vec{u} \in \vec{H}^{\vec{s}}$.   Then $ \overline{\vec{u}} \in
\vec{H}^{{s}+ \beta}$  and 
\begin{equation}
\label{1suralpharbeta}
\|\overline{\vec{u}}  \|_{\vec{H}^{{s}+ \beta}}^r  \le \frac{1}{\alpha^{r \beta}} \|\vec{u} \|_{\vec{H}^{{s}}}^r.
  \end{equation}
\end{lemme}
\textbf{Proof.} 
When $ \displaystyle \bu=\sum_{{\vec{k}}\in \vec{I}}\hat{\vec{u}}_{\vec{k}}\exp\left\lbrace   i\vec{k}\cdot \vec{x}\right\rbrace $, then 
\BEQ \label {GHGGGPRIME}  \overline{\bu} =\sum_{{\vec{k}}\in {\vec{I}}}  \frac{\hat{\bu}_{\vec{k}}}{1 + \alpha^{2 \theta} | {\vec{k}}|^{2\theta}}    \exp\left\lbrace   i\vec{k}\cdot \vec{x}\right\rbrace.
\EEQ
Formula (\ref{GHGGGPRIME}) easily yields   the estimate 
\begin{equation}
\begin{array}{ccccccc}
\label{1suralpharbetaprime}
\displaystyle \|\overline{\vec{u}}  \|_{\vec{H}^{{s}+ \beta}}^r  &\le& \displaystyle \left(\sup_{\vec{k} \in \vec{I}}{\frac{|\vec{k}|^{2 \beta}}{(1+\alpha^{2\theta}|\vec{k}|^{2 \theta})^{2}}}\right)^{\frac{r}{2}} \|\vec{u} \|_{\vec{H}^{{s}}}^r\\
\displaystyle &\le& \displaystyle \left(\frac{1}{\alpha^{2 \beta}}\sup_{\vec{k} \in \vec{I}}\left({\frac{\alpha^{\beta}|\vec{k}|^{\beta}}{1+\alpha^{2\theta}|\vec{k}|^{2\theta}}}\right)^2\right)^{\frac{r}{2}} \|\vec{u} \|_{\vec{H}^{{s}}}^r\\
 &\le& \displaystyle \frac{1}{\alpha^{r \beta}} \|\vec{u} \|_{\vec{H}^{{s}}}^r,\\
\end{array}
  \end{equation}
  where we have used the fact that 
  \begin{equation}
\sup_{\vec{k} \in \vec{I}}\left({\frac{\alpha^{\beta}|\vec{k}|^{\beta}}{1+\alpha^{2\theta}|\vec{k}|^{2\theta}}}\right)^2 \le 1 , \ \hbox{for all }  0 \le \beta \le 2\theta.
  \end{equation}
  
\begin{lemme}
Assume  $\vec{u}_{\alpha}$ belongs to $  $ the enegy space of  solutions of the Navier-Stokes equations, then 
\begin{align}
\overline{\bu_\alpha}  &\rightarrow \bu &&\textrm{strongly in  }
L^p(0,T;L^p(\tore)^3) \textrm{ for all  }   2 \le p < \frac{10}{3}.
\end{align}
\end{lemme}
\textbf{Proof.} 
From the Sobolev injection $ \displaystyle \vec{H}^{\frac{3p-6}{2p}}\hookrightarrow L^p(\tore)^3,$ it is sufficent to show that 
  \begin{equation}
\int_{0}^{T}\| \overline{\bu_\alpha} -  \bu_{\alpha}\|_{\vec{H}^{\frac{3p-6}{2p}}}^p dt \longrightarrow 0, \hbox{ when } \alpha \rightarrow 0. 
  \end{equation}
  From the relation between $\overline{\bu_\alpha} $  and $\bu_{\alpha} $  we have 
    \begin{equation}
  \| \overline{\bu_\alpha} -  \bu_{\alpha}\|_{\vec{H}^{\frac{3p-6}{2p}}}^p \le \alpha^{2 \theta p} \|\overline{\bu_\alpha} \|_{\vec{H}^{\frac{3p-6}{2p} + 2\theta}}^p
   \end{equation}
Lemma \ref{fourierdiscret} implies that 
 \begin{equation}
\int_{0}^{T}\| \overline{\bu_\alpha} -  \bu_{\alpha}\|_{\vec{H}^{\frac{3p-6}{2p}}}^p dt \le \alpha^{\frac{10-3p}{2}} \int_{0}^{T}\|{\bu_\alpha} \|_{\vec{H}^{\frac{2}{p}}}^p dt.
 \end{equation}
Recall that   \begin{equation}
 \int_{0}^{T}\|{\bu_\alpha} \|_{\vec{H}^{\frac{2}{p}}}^p dt < \infty \hbox{ fo any } 2 \le p \le \infty.  
   \end{equation}
This yields  the desired result for any $ 2 \le p  < 10/3$.


In order to show that $(\vec{u}_{\alpha}, p_{\alpha})$ gives rise to a   suitable solution of the  Navier-Stokes equations, it is necessary to take the   limit  $\alpha \rightarrow 0$. 
We have the following theorem:
\begin{Theorem}
Let $(\vec{u}_{\alpha}, p_{\alpha})$ be the solution of (\ref{alpha ns}), in the $ \alpha \rightarrow 0$ limit, $(\vec{u}_{\alpha}, p_{\alpha})$ tends to a suitable solution to the Navier-Stokes equations. 
\end{Theorem}
\textbf{Proof:} 
 By a classical compactness argument \cite{FDT02}, 
we deduce  that $\vec{u}_{\alpha}$ approaches $\vec{u}$ strongly in $L^2([0,T],\vec{H})$ for all $  T >0$ when $\alpha$ tends to zero,
where $\vec{u}$ is a weak solution to the Navier-Stokes equations . 
  It remains to show that $(\vec{u}_{}, p_{})$  verifies the following local energy inequality   
\begin{equation}
\label{local}
\begin{array}{cccc}
\displaystyle 2\nu\int_{0}^{T}\int_{\tore}|\nabla \vec{u}_{}|^{2}\phi \ d\vec{x}dt \le \displaystyle \int_{0}^{T}\int_{\tore}| \vec{u}_{}|^{2}\left(\phi_t + \nu \Delta \phi \right) + \left(| \vec{u}_{}|^{2} {\vec{u}_{}}+ 2p\vec{u}_{}\right) \cdot \nabla \phi   \  d\vec{x}dt\\
 \hskip 4cm \displaystyle + \displaystyle 2\int_{0}^{T}\int_{\tore} \vec{f} \vec{u}\phi \  d\vec{x}dt. 
 \end{array}
\end{equation}
To do so, we need to find estimates
that are independent from $\alpha$. 
Using the fact that  $(\vec{u}_{\alpha}, p_{\alpha})$  belong to the energy space:  $L^{\infty}
([0,T];\vec{H})\cap L^2([0,T];\vec{V}) \cap L^{\frac{10}{3}}([0,T];L^{\frac{10}
{3}}(\tore)^3) $ and from the  Aubin-Lions compactness Lemma  (the same arguments as in section 4) we can find a not relabeled subsequence $(\vec{u}_{\alpha}, p_{\alpha})$  and $(\vec{u}_{}, p_{})$ such that when   $\alpha$ tends to zero we have: 
\begin{align}
\bu_\alpha &\rightharpoonup^* \bu &&\textrm{weakly$^*$ in } L^{\infty}
([0,T];\vec{H}), \label{c122}\\
\bu_\alpha &\rightharpoonup \bu &&\textrm{weakly in }
L^2([0,T];\vec{V})\cap L^{\frac{10}{3}}([0,T];L^{\frac{10}
{3}}(\tore)^3), \label{c22}\\
p_\alpha&\rightharpoonup p &&\textrm{weakly in } L^{\frac{5}{3}}([0,T];L^{\frac{5}{3}}
(\tore)), \label{c32}\\
\bu_\alpha &\rightarrow \bu &&\textrm{strongly in  }
L^p([0,T];L^p(\tore)^3) \textrm{ for all } 2 \le  p<\frac{10}{3},\label{c82}\\
\overline{\bu_\alpha}  &\rightarrow \bu &&\textrm{strongly in  }
L^p([0,T];L^p(\tore)^3) \textrm{ for all } 2 \le p < \frac{10}{3},\\
p_\alpha&\rightharpoonup p &&\textrm{strongly in  }
L^p([0,T];L^p(\tore)) \textrm{ for all } 1 < p < \frac{5}{3}.
\end{align}

These convergence results allow us  to take the limit in (\ref{local alpha}) in order to obtain, by using the weak lower semicontinuty of the norm in $L^2([0,T];\vec{V})$
$$ \liminf_{\alpha\rightarrow0}\int_{0}^{T}\|\nabla \vec{u}_{\alpha}\|_{\vec{H}}^2 \ dt  \ge \int_{0}^{T}\|\nabla \vec{u}\|_{\vec{H}}^2 \ dt,$$
that  $(\bu,p)$  verifies  the local energy inequality (\ref{local}).

%

\section{The deconvolution case}

\subsection{The modified deconvolution operator}\label{sec:dec}
In this section, we introduce the modified deconvolution operator which interpolate the usual  deconvolution operator introduced in \cite{LL08}.

Let $\alpha >0$, $s \ge -1$, $  0 \le \theta  \le 1$, $\bu \in \vec{H}^{s}$ and let  $\overline{\bu} \in \vec{H}^{s+2\theta}$ be the unique solution to the
 equations 
\BEQ \label {VPYRA} \alpha^{2 \theta}(-\Delta)^{\theta} \overline{\bu} + \overline{\bu} = \bu. \EEQ
\BEQ
\nabla \cdot \vec{u}= \nabla \cdot \overline{\vec{u}}=0
\EEQ
We also shall  denote by $\mathbb{G}$ the operator 
\BEQ \mathbb{G} : \begin{array} {l} \vec{H}^{s+2\theta} \longrightarrow {\vec{H}^{s}}, \\
\bw \longrightarrow   \alpha^{2 \theta}(-\Delta)^{\theta} \bw  + \bw. \end{array}  \EEQ
Therefore, one has 
\BEQ \overline{\bu} = \mathbb{G} ^{-1} \bu. \EEQ


Let us consider the operators 
$$ D_N = \sum_{n = 0}^N ( I - \mathbb{G} ^{-1} )^n.$$
and 
\BEQ \label{hn} {H_N} (\bu) = D_N (\overline{\bu}).  \EEQ
A straightforward calculation yields 
 \BEQ H_N \left ( \sum_{ {\bf k}\in {\cal  I}_3} \bu_{\bf k} \expk \right ) =  
\sum_{ {\bf k}\in {\cal  I}_3} \left ( 1- \left ( \frac{\alpha^{2 \theta} | {\bf k} |^{2 \theta}}
 {1 +  \alpha^{2 \theta} | {\bf k} |^{2 \theta}}\right )^{N+1} \right )  \bu_{\bf k} \expk. \EEQ 
One can prove the following (see in \cite{LL08, RL} with $ \theta=1$): 
\begin{itemize} 
\item Assume $\bu \in {\vec{H}^{s}}$. Then $H_N (\bu) \in \vec{H}^{s+2\theta}$ and 
$ || H_N (\bu) ||_{\vec{H}^{s+2\theta}}\le C(N, \alpha) || \bu ||_{\vec{H}^{s}}$, where $C(N, \alpha)$ blows up when $\alpha$ goes to zero and/or $N$ goes to infinity. This is due to the fact 
$$ \left ( 1- \left ( \frac{\alpha^{2 \theta} | {\bf k} |^{2 \theta}}{1 +  \alpha^{2 \theta} | {\bf k} |^{2 \theta}} \right )^{N+1} \right )  \approx \frac{N+1}{\alpha^{2 \theta} | {\bf k} |^{2 \theta} }\quad \hbox{as} \quad | {\bf k} |_\infty \rightarrow \infty. $$

\item The operator $H_N$ maps continuously ${\vec{H}^{s}}$ into ${\vec{H}^{s}}$ and 
$ || H_N ||_{ {\cal L} ({{H}^{s},H^s}) } = 1$. 

\item Assume $\bu \in {\vec{H}^{s}}$. Then the sequence $(\overline{\vec{u}}_\alpha)_{\alpha >0}$ converges strongly to $\bu$ in the space ${\vec{H}^{s}}$.

\item Let $\bu \in {\vec{H}^{s}}$. Then the sequence $(H_N (\bu) )_{N \in \N}$ converges strongly to $\bu$ in ${\vec{H}^{s}}$ when $N$ goes to infinity. 
\end{itemize} 
\subsection{The deconvolution model: well-posedness results}
\label{sec:deconvolution}
We consider here a family of equations interpolating between the Navier-Stokes equations \cite{JL34} and the Leray-deconvolution model \cite{LL08}  with periodic boundary conditions.

\begin{equation}
\label{dalpha ns}
 \left\{
 \begin{array} {llll} \displaystyle
 \frac{\partial \vec{u}^{}}{\partial t}+ H_{N}(\vec{u})\cdot \nabla
\vec{u}^{} - \nu \Delta \vec{u}^{} + \nabla
p^{} = \vec{f} \ \ \ \ \ \hbox{in}\ \RR^{+}\times\mathbb{T}_3,\\
\nabla \cdot \vec{u}^{}= \nabla \cdot H_{N}({\vec{u}})=0,\
 \displaystyle \int_{\mathbb{T}_3} \vec{u}^{}=\displaystyle \int_{\mathbb{T}_3}H_{N}({\vec{u}})=0,\\
\vec{u}^{}(t,\vec{x}+{L}\vec{e}_{\vec{j}})=\vec{u}^{}(t,\vec{x}),\\
\vec{u}^{}_{t=0}=\vec{u}_{0}^{}.
\end{array}
\right.
\end{equation}

Where $H_{N}(\vec{u})$ is an ineterpolating deconvolution operater introduced above.\\

When $N=0$,we obtain the equations (\ref{alpha ns}).
Similarly, we obtain  the same results of well-posedness for the interpolating model (\ref{dalpha ns}). 

\begin{Theorem}
\label{2hn}
  For any $T>0$,  let $\vec{f} \in L^{2}([0,T],\vec{V}^{-{1}})$ and $\vec{u}_{0} \in \vec{H}$. Assume that  $\theta= 1/4 $. Then there exists a unique     solution $(\vec{u}_{},p_{}):=(\vec{u}_{\alpha,N},p_{\alpha,N})$ to eqs. (\ref{dalpha ns}) that satisfies 
$ \vec{u}_{} \in  C([0,T];\vec{H}) \cap L^{2}([0,T];\vec{V})$ and $\displaystyle \frac{d\vec{u}}{dt} \in L^{2}([0,T];\vec{V^{-1}})$ and  ${p}_{} \in  L^{2}([0,T],L^{{2}}(\mathbb{T}_3))$.\\
Such that $\vec{u}_{}$ verifies 
$$\left\langle  \frac{d \vec{u}_{}}{d t}  +\nu A \vec{u}_{} + B(H_N({\vec{u}}),\vec{u})
 -\vec{f}, \phi  \right\rangle_{\vec{V}^{-{{1}{}}},\vec{V}^{{{1}{}}}}=0 $$
 for every $\phi \in \vec{V}^{}$ and almost every $ t \in (0,T)$.
Moreover, this solution depends continuously on the initial data $\vec{u}_{0}$ in the $L^{{2}}$ norm.
\end{Theorem}
And 
\begin{Theorem}
\label{1hn}
For any $T>0$, let $\vec{f} \in L^{2}([0,T],\vec{V}^{-{\frac{1}{2}}})$ and $\vec{u}_{0} \in \vec{V}^{\frac{1}{2}}$. Assume that $\theta= \frac{1}{4}$. Then  there exists a unique solution  $(\vec{u}_{},p_{}):=(\vec{u}_{\alpha,N},p_{\alpha,N})$ to eqs. (\ref{dalpha ns})  that satisfies 
 $ \vec{u}_{} \in  C([0,T];\vec{V}^{\frac{1}{2}}) \cap L^{2}([0,T];\vec{V}^{\frac{3}{2}})$, $\displaystyle \frac{d\vec{u}_{}^{}}{dt} \in L^{2}([0,T],\vec{V}^{{-\frac{1}{2}}})$ and  ${p}_{} \in  L^{2}([0,T],H^{\frac{1}{2}}(\mathbb{T}_3))$.\\
Such that  $\vec{u}$ satisfies 
$$\left\langle  \frac{d \vec{u}^{}}{d t}  +\nu A \vec{u}^{} + B(H_N({\vec{u}}),\vec{u})
 -\vec{f}, \phi  \right\rangle_{\vec{V}^{-\frac{1}{2}},\vec{V}^{\frac{1}{2}}}=0 $$
 for every $\phi \in \vec{V}^{\frac{1}{2}}$ and almost every $ t \in (0,T)$.
Moreover, this solution depends continuously on the initial data $\vec{u}_{0}$ in the $H^{\frac{1}{2}}$ norm.
\end{Theorem}
Since the proof of Theorem \ref{2hn} and \ref{1hn} are similar to that of Theorem \ref{2} and \ref{1} we omit their proof here.\\

\subsection{Convergence to a suitable weak solution to the NSE}
By a similar argument to that in \cite{LL08}, it is possible to prove that up to a subsequence the solution for  the deconvolution model (with $\theta=\frac{1}{4}$) converges when $\alpha \rightarrow 0$ or/and $N\rightarrow +\infty$ to a weak Leray solution to the  Navier-Stokes equations. Next, we will show that the deconvolution regularization  gives rise to a suitable weak solution to the Navier-Stokes equations. We need first 
\begin{lemme}
Let $(\vec{u}_{\alpha,N}, p_{\alpha,N})$ be the unique solution of  (\ref{dalpha ns}) with  $\theta \ge \frac{1}{4}$ then  $(\vec{u}_{\alpha,N}, p_{\alpha,N})$ verifies the following local inequality
\begin{equation}
\label{dlocal alpha}
\begin{array}{ccc}
\displaystyle 2\nu\int_{0}^{T}\int_{\tore}|\nabla \vec{u}_{\alpha,N}|^{2}\phi \ d\vec{x}dt = \displaystyle \int_{0}^{T}\int_{\tore}| \vec{u}_{\alpha,N}|^{2}\left(\phi_t + \nu \Delta \phi \right)\  d\vec{x}dt \\ +\displaystyle \int_{0}^{T}\int_{\tore} \left(| \vec{u}_{\alpha,N}|^{2} H_N({\vec{u}_{\alpha,N}})+ 2p_{\alpha,N}\vec{u}_{\alpha,N}\right) \cdot \nabla \phi   \  d\vec{x}dt
\displaystyle + \displaystyle 2\int_{0}^{T}\int_{\tore} \vec{f} \vec{u}_{\alpha,N}\phi \  d\vec{x}dt, 
\end{array}
\end{equation}
for all $T \in (0,+\infty]$ and for all  non negative fonction  $\phi \in C^{\infty}$ and  supp $\phi \subset\subset \tore\times(0,T).$
\end{lemme}
\textbf{Proof.} See Lemma \ref{lemmelocalen}.\\

Then we have 
\begin{lemme}
\label{conlp}
Assume  $\vec{u}_{\alpha,N}$ belong to $  $ the enegy space of the solutions of the Navier-Stokes equations, then 
\begin{align}
H_N({\bu_{\alpha,N}})  &\rightarrow \bu &&\textrm{strongly in  }
L^p(0,T;L^p(\tore)^3) \textrm{ for all  }   2 \le p < \frac{10}{3}.
\end{align}
\end{lemme}
\textbf{Proof.} 
From the sobolev injection $ \displaystyle \vec{H}^{\frac{3p-6}{2p}}\hookrightarrow L^p(\tore)^3,$ it is sufficent to show that 
  \begin{equation}
\int_{0}^{T}\| H_N({\bu_{\alpha,N}}) -  \bu_{\alpha,N}\|_{\vec{H}^{\frac{3p-6}{2p}}}^p dt \longrightarrow 0, \hbox{ when } \alpha \rightarrow 0 \hbox{ or/and } N \rightarrow \infty. 
  \end{equation}
  From the relation between $H_N(\bu_{\alpha,N}) $  and $\bu_{\alpha,N} $  we have 
    \begin{equation}
    \begin{array}{ccccc}
  \displaystyle \| H_N({\bu_{\alpha,N}}) -  \bu_{\alpha,N}\|_{\vec{H}^{\frac{3p-6}{2p}}}^p = \displaystyle \sum_{\vec{k} \in \vec{I}}\left(|\vec{k}|^{\frac{3p-6}{p}}\left(  \frac{\alpha^{2\theta} \vec{k}^{2\theta}}{1+ \alpha^{2\theta} \vec{k}^{2\theta}}\right)^{2(N+1)}|\hat{\vec{u}}_{\vec{k}}|^2\right)^{\frac{p}{2}}\\
 \hskip 4cm \le  \displaystyle \left(\sup_{\vec{k} \in \vec{I}}\left(  |\vec{k}|^{\frac{3p-10}{p}}  \left(\frac{\alpha^{2\theta}|\vec{k}|^{2\theta}}{1+\alpha^{2\theta}|\vec{k}|^{2\theta}}\right)^{2(N+1)}\right)\right)^{\frac{p}{2}} \| \vec{u}\|_{\vec{H}^{\frac{2}{p}}}^p
  \end{array}
   \end{equation}
   Using the fact that $\displaystyle l^m \hookrightarrow l^{\infty}$ for any $m \ge 1$ we get for $p < \frac{10}{3}$ that 
      \begin{equation}
    \begin{array}{ccccc}
  \displaystyle \| H_N({\bu_{\alpha,N}}) -  \bu_{\alpha,N}\|_{\vec{H}^{\frac{3p-6}{2p}}}^p 
 \le  \displaystyle \left(\sum_{\vec{k} \in \vec{I}}\left( \frac{1}{ |\vec{k}|^{m\frac{10-3p}{p}}}  \left(\frac{\alpha^{2\theta}|\vec{k}|^{2\theta}}{1+\alpha^{2\theta}|\vec{k}|^{2\theta}}\right)^{2m(N+1)}\right)\right)^{\frac{p}{2}} \| \vec{u}\|_{\vec{H}^{\frac{2}{p}}}^p
  \end{array}
   \end{equation}
   Now, we choose $m$ such that $m\frac{10-3p}{p} > 3$, and recall  that   \begin{equation}
 \int_{0}^{T}\|{\bu_{\alpha,N}} \|_{\vec{H}^{\frac{2}{p}}}^p dt < \infty \hbox{ fo any } 2 \le p \le \infty,  
   \end{equation}
 and 
  \begin{equation}
  \frac{\alpha^{2\theta}|\vec{k}|^{2\theta}}{1+\alpha^{2\theta}|\vec{k}|^{2\theta}}  \displaystyle \stackrel{\alpha \rightarrow 0}{\longrightarrow 0}  \hbox{ or/and } 
\sup_{\vec{k} \in \vec{I}}\left({\frac{\alpha^{2\theta}|\vec{k}|^{2\theta}}{1+\alpha^{2\theta}|\vec{k}|^{2\theta}}}\right)^{2m} \le 1.
  \end{equation}
  Using the dominate convergence theorem, we conclude
 the desired result for any $ 2 \le p  < 10/3$.\\
 \begin{Remark}
 Note that the above result holds true for $p=10/3$.
 \end{Remark}
  The above $L^p$ convergence combined with the fact that  $\vec{u}_{\alpha,N}$ belong to  the energy space of solutions of the Navier-Stokes equations and the  Aubin-Lions compactness  Lemma allow us to take the limit $ \alpha \rightarrow 0$ or/and $ N \rightarrow  \infty$ in (\ref{dlocal alpha})  and to deduce the following Theorem:
  \begin{Theorem}
Let  $(\vec{u}_{\alpha,N}, p_{\alpha,N})$ be the solution of  (\ref{dalpha ns}), Then when  $\alpha \rightarrow 0$ and/or $N \rightarrow \infty$,  $(\vec{u}_{\alpha,N}, p_{\alpha, N})$ tends to a suitable weak solution  $(\vec{u}_{}, p_{})$ to the  Navier-Stokes equations. The convergence to $\vec{u}$ is weak in $L^2([0,T];\vec{V})\cap L^{\frac{10}{3}}([0,T];L^{\frac{10}
{3}}(\tore)^3)$ and strong in \\
$L^q([0,T];L^q(\tore)^3)$  for all $ 2 \le  q<\frac{10}{3}$ and the convergence to $p$ is strong  in $L^q([0,T];L^q(\tore))$  for all $ 1 <  q<\frac{5}{3}$ and weak  in  $L^{\frac{5}{3}}([0,T];L^{\frac{5}{3}}(\tore))$.
\end{Theorem}

\section{The MHD Case}\label{sec:MHD}
In this section, we consider a deconvolution-type regularization of the  magneto-hydrodynamic (MHD) equations, given by
\begin{subequations}\label{IMHDV}
\begin{align}
 \partial_t\bu-\nu_1 \Delta \bu +(H_N(\bu)\cdot\nabla)\bu-(H_N(\B)\cdot\nabla)\B
+\nabla p+\frac{1}{2}\nabla|\B|^2&=0,
\label{IMHDV1}\\
\partial_t\B -\nu_2 \triangle \B +(H_N(\bu)\cdot\nabla)\B-(H_N(\B)\cdot\nabla)\bu&=0,
\label{IMHDV2}\\
 \nabla\cdot \B = \nabla\cdot H_N(\B) =\nabla\cdot \bu= \nabla\cdot H_N(\bu)&= 0,\label{IMHDV3}\\
 \B(0)=\B_{0},\; \bu(0)=\bu_{0},&\label{IMHDV4}
\end{align}
\end{subequations}
where the boundary conditions are taken to be periodic, and we also assume as before that
$\displaystyle \int_{\tore} \bu\,dx =  \int_{\tore} \B\,dx =0.$\\
 Here, the unknowns are the fluid velocity field $\bu(t,\vec{x})$, the fluid pressure $p(t,\vec{x})$, and the magnetic field $\B(t,\vec{x})$.  Note that when $\alpha=0$, we formally retrieve the  MHD equations and the  MHD-Leray-$\alpha$ equations  are obtained when $N=0$ and $\theta=1$.\\
Existence and uniqueness results for MHD equations are established by G. Duvaut and J.L. Lions in \cite{DuLi72}. These results are completed by M. Sermange and R. Temam in \cite{RT83mhd}. They showed that the classical properities of the Navier-Stokes equations can be extended to the MHD system. 
 The use of Leray-$\alpha$ regularization to the MHD equations has received many studies see \cite{LT2007}. The idea to use the deconvolution operator from \cite{LL08} to MHD equations is a new feature for the present work.\\
\subsection{Existence, unicity and convergence results}
First, we establish  the global existence and uniqueness of solutions for the MHD-Deconvolution equations (\ref{IMHDV}) for $\theta=1/4$.\\
We have the following theorem:
\begin{Theorem}
\label{existensssIMHDV}
 For $\theta=1/4$.  Assume ${\vec{u}}^{}_{0} \in \vec{H} $ and $\B_{0} \in \vec{H}
$.  Then for any $T>0$, (\ref{IMHDV}) has a unique regular    solution $ (\vec{u}^{ },\B, {p}^{ }):=(\vec{u}_{\alpha,N },\B_{\alpha, N},{p}_{\alpha,N })$
 such that,
 ${\vec{u}^{ }}, {\B^{}} \in C((0,T],\vec{H})\cap L^{2}([0,T],\vec{V}^{1}),$ and  
 ${p}_{}  \in  L^{2}([0,T],L^{{2}}(\mathbb{T}_3)).$\\
 Furthermore  the solution verifies  
 \begin{equation}
\label{localMHD}
\begin{array}{cccc}
\hskip -3cm \displaystyle 2\nu\int_{0}^{T}\int_{\tore}\left(|\nabla \vec{u}_{}|^{2}+|\nabla \B_{}|^{2} \right)\phi \ d\vec{x}dt \\
\hskip 1cm = \displaystyle \int_{0}^{T}\int_{\tore}| \vec{u}_{}|^{2}\left(\phi_t + \nu_1 \Delta \phi \right) +| \B_{}|^{2}\left(\phi_t + \nu_2 \Delta \phi \right) \  d\vec{x}dt\\
\hskip 1.5cm +\displaystyle\int_{0}^{T}\int_{\tore}  \left( \left(| \vec{u}_{}|^{2}+| \B_{}|^{2}\right) {H_N(\vec{u}_{}})+ 2p\vec{u}_{}\right) \cdot \nabla \phi   \  d\vec{x}dt\\
 \hskip 2cm 
 \displaystyle -\displaystyle 4\int_{0}^{T}\int_{\tore}  \left( \vec{u}_{} \B_{}\right){H_N(\B_{}})\cdot \nabla \phi\  d\vec{x}dt 
 \end{array}
\end{equation}
for all $T \in (0,+\infty]$ and for all  non negative fonction  $\phi \in C^{\infty}$ and  supp $\phi \subset\subset \tore\times(0,T).$
\end{Theorem}
 \textbf{Proof.} 
 We only sketch  the proof since  is similar to the   Navier-Stokes equations case. The proof  is obtained by taking the inner product of (\ref{IMHDV1}) with $\vec{u},$ (\ref{IMHDV2}) with $\B$ and then adding them, the existence of a
solution to Problem $(\ref{IMHDV})$ can be derived thanks to the Galerkin
method. Notice that  ${\vec{u}^{}}, {\B^{ }}$    satisfy the following estimates
\begin{equation}
  \displaystyle   \frac{1}{2} \frac{d}{dt}\left(
\|{\vec{u}}^{ }(t,\vec{x})\|_{\vec{H}}^{2}  +
\|{\B}^{ }(t,\vec{x})\|_{\vec{H}}^{2}\right) +\min{(\nu_1,\nu_2)} \left(\|\nabla
{\vec{u}}^{ }(t,\vec{x})\|_{\vec{H}}^{2} +\|\nabla
{\B}^{ }(t,\vec{x})\|_{\vec{H}}^{2}\right) 
\le 0
\end{equation}
The pressure $p^{ }$   is reconstructed from $\vec{u}, H_N(\bu), \B, H_N(\B)$  (as we work with periodic boundary conditions) 
and its regularity
results from the fact that  $(H_N(\bu)\cdot\nabla)\bu, (H_N(\B)\cdot\nabla)\B \in L^2([0,T],\vec{H}^{-1})$.\\
It remains to prove the uniqueness. Let $(\vec{u}_{1},\B_{1}, p_{1})$ and $(\vec{u}_{2},\B_{2},
p_{2})$,  be two solutions, $\delta\vec{u} = \vec{u}_{2}-\vec{u}_{1}$,$\delta\B = \B_{2}-\B_{1}$, $\delta p =
p_{2}-p_{1}$. Then one has
\begin{subequations}
\begin{align}
\label{matin1}
\partial_{t} \delta\vec{u} + (H_{N} (\vec{u}_{1}) \nabla) \delta\vec{u}-(H_{N} (\B_{1}) \nabla) \delta\B  - \nu_1
\Delta\delta\vec{u} + \nabla\delta p = -(H_{N} (\delta\vec{u}) \nabla) \vec{u}_{2} +(H_{N} (\delta\B) \nabla) \B_{2},\\
\label{matin2}
\partial_{t} \delta\B +(H_{N} (\vec{u}_{1}) \nabla) \delta\B -  (H_{N} (\B_{1}) \nabla) \delta\vec{u} - \nu_2
\Delta\delta\B  = (H_{N} (\delta\B) \nabla) \vec{u}_{2}- (H_{N} (\delta\vec{u}) \nabla) \B_{2},
\end{align}
\end{subequations}
and $\delta\vec{u} = 0$, $\delta\B =0$  at initial time. One can take $\delta\vec{u} \in  L_{}^{\infty}((0,T],\vec{H})\cap L^{2}([0,T],\vec{V})$
 as test in (\ref{matin1}) and $\delta\B \in  L_{}^{\infty}((0,T],\vec{H})\cap L^{2}([0,T],\vec{V})$
 as test in (\ref{matin2}). Since $H_{N}
(\vec{u}_1)$ is divergence-free, one has
 \begin{subequations}
\begin{align}\displaystyle
\int_{0}^{T} \int_{\tore}(H_{N} (\vec{u}_{1}) \nabla) \delta\vec{u}. \delta\vec{u} =
0,\\
\displaystyle
\int_{0}^{T} \int_{\tore}(H_{N} (\vec{u}_{1}) \nabla) \delta\B. \delta\B =
0.
\end{align}
\end{subequations}

Therefore, 

 \begin{equation}
\label{soira}\begin{array}{ll}
 \displaystyle {\frac{d }{2dt}} \int_{\tore} | \delta\vec{u}| ^{2} +\displaystyle \nu_1\int_{\tore}| \nabla
\delta\vec{u} |^{2}  -\displaystyle \int_{\tore} (H_{N} (\B_{1}) \nabla) \delta\B. \delta
\vec{u}\\
 \hskip 3cm  =\displaystyle -\int_{\tore}(H_{N} (\delta\vec{u}) \nabla) \vec{u}_{2}. \delta
\vec{u} +\displaystyle \displaystyle\int_{\tore}(H_{N} (\delta\B) \nabla) \B_{2}. \delta
\vec{u},
\end{array}
\end{equation}
and 
  \begin{equation}
\label{soirb}\begin{array}{ll}
 \displaystyle {\frac{d }{2dt}} \int_{\tore}  | \delta\B| ^{2} + \displaystyle\nu_2\int_{\tore} | \nabla
\delta\B |^{2}   - \displaystyle \int_{\tore} (H_{N} (\B_{1}) \nabla) \delta\vec{u} . \delta
\B\\ 
 \hskip 3cm  = \displaystyle \int_{\tore}(H_{N} (\delta\B) \nabla) \vec{u}_{2}     . \delta
\B   \displaystyle   -\displaystyle\int_{\tore}   (H_{N} (\delta\vec{u}) \nabla) \B_{2}  . \delta
\B.
\end{array}
\end{equation}

One has by a integration by parts,
\BEQ \label{MATIN2}
- \displaystyle \int_{\tore} (H_{N} (\B_{1}) \nabla) \delta\B. \delta
\vec{u}     = \displaystyle \int_{\tore} (H_{N} (\B_{1}) \nabla) \delta\vec{u} . \delta
\B. 
\EEQ
One has by adding (\ref{soira}),(\ref{soirb}) and using (\ref{MATIN2})

 \begin{equation}
 \label{soirab}
\begin{array}{ll}
 \displaystyle {\frac{d }{2dt}} \int_{\tore} | \delta\vec{u}| ^{2}+{\frac{d }{2dt}} \int_{\tore}  | \delta\B| ^{2} +\displaystyle \nu_1\int_{\tore}| \nabla
\delta\vec{u} |^{2} +\displaystyle\nu_2\int_{\tore} | \nabla
\delta\B |^{2} \\
\hskip 2cm=\displaystyle -\int_{\tore}(H_{N} (\delta\vec{u}) \nabla) \vec{u}_{2}. \delta
\vec{u} +\displaystyle \displaystyle\int_{\tore}(H_{N} (\delta\B) \nabla) \B_{2}. \delta \vec{u} \\
\hskip 4cm
+\displaystyle \int_{\tore}(H_{N} (\delta\B) \nabla) \vec{u}_{2}     . \delta
\B   \displaystyle   -\displaystyle\int_{\tore}   (H_{N} (\delta\vec{u}) \nabla) \B_{2}  . \delta
\B.
 \end{array}
\end{equation}

One has by a  integration by parts,
\begin{subequations}
\begin{align}
-\int_{\tore}(H_{N} (\delta\vec{u}) \nabla) \vec{u}_{2}. \delta
\vec{u} =\int_{\tore} H_{N} (\delta\vec{u}) \otimes \vec{u}_{2} : \g \delta
\vec{u},\\
\int_{\tore}(H_{N} (\delta\B) \nabla) \B_{2}. \delta
\vec{u} = -\int_{\tore} H_{N} (\delta\B) \otimes \B_{2} : \g \delta
\vec{u},\\
\int_{\tore}(H_{N} (\delta\B) \nabla) \vec{u}_{2}. \delta
\B = -\int_{\tore} H_{N} (\delta\B) \otimes \vec{u}_{2} : \g \delta
\B,\\
-\int_{\tore}(H_{N} (\delta\vec{u}) \nabla) \B_{2}. \delta
\B = \int_{\tore} H_{N} (\delta\vec{u}) \otimes \B_{2} : \g \delta
\B.
\end{align}
\end{subequations}

By Young's inequality, 

\begin{subequations}
\begin{align}
|\int_{\tore}(H_{N} (\delta\vec{u}) \nabla) \vec{u}_{2}. \delta
\vec{u}|  \le \frac{\nu_1}{4} \int_{\tore} |\g \delta \vec{u} |^2 + \frac{1}{\nu_1} \int _{\tore}| H_N (\delta \vec{u}) |^2 | \vec{u}_2|^2,\\
 |\int_{\tore}(H_{N} (\delta\B) \nabla) \B_{2}. \delta
\vec{u}|  \le \frac{\nu_1}{4} \int_{\tore} |\g \delta \vec{u} |^2 + \frac{1}{\nu_1} \int _{\tore}| H_N (\delta \B) |^2 | \B_2|^2,\\
|\int_{\tore}(H_{N} (\delta\B) \nabla) \vec{u}_{2}. \delta
\B|  \le \frac{\nu_2}{4} \int_{\tore} |\g \delta \B |^2 + \frac{1}{\nu_2} \int _{\tore}| H_N (\delta \B) |^2 | \vec{u}_2|^2,\\
|\int_{\tore}(H_{N} (\delta\vec{u}) \nabla) \B_{2}. \delta
\B|  \le \frac{\nu_2}{4} \int_{\tore} |\g \delta \B |^2 + \frac{1}{\nu_2} \int _{\tore}| H_N (\delta \vec{u}) |^2 | \B_2|^2.
 \end{align}
\end{subequations}

By H\"{o}lder inequality combined with Sobolev injection
\begin{subequations}
\begin{align}
\begin{array}{llll}
 \displaystyle \frac{1}{\nu_1} \int _{\tore}| H_N (\delta \vec{u}) |^2 | \vec{u}_2|^2 &\le \displaystyle \frac{1}{\nu_1} \| H_N (\delta \vec{u}) \|_{L^3}^2 \| \vec{u}_2\|_{L^6}^2&\\
 & \displaystyle\le  \frac{1}{\nu_1} \| H_N (\delta \vec{u}) \|_{\vec{V}^{\frac{1}{2}}}^2 \| \vec{u}_2\|_{\vec{V}}^2,&
\end{array} \\
\begin{array}{llll}
 \displaystyle \frac{1}{\nu_1} \int _{\tore}| H_N (\delta \B) |^2 | \B_2|^2
 &\le \displaystyle \frac{1}{\nu_1} \| H_N (\delta \B) \|_{L^3}^2 \| \B_2\|_{L^6}^2&\\
 & \displaystyle\le  \frac{1} {\nu_1} \| H_N (\delta \B) \|_{\vec{V}^{\frac{1}{2}}}^2 \| \B_2\|_{\vec{V}}^2,&
 \end{array}\\
\begin{array}{llll}
 \displaystyle \frac{1}{\nu_2} \int _{\tore}| H_N (\delta \B) |^2 | \vec{u}_2|^2  
 &\le \displaystyle \frac{1}{\nu_2} \| H_N (\delta \B) \|_{L^3}^2 \| \vec{u}_2\|_{L^6}^2&\\
 & \displaystyle\le  \frac{1}{\nu_2} \| H_N (\delta \B) \|_{\vec{V}^{\frac{1}{2}}}^2 \| \vec{u}_2\|_{\vec{V}}^2,&
 \end{array}\\
\begin{array}{llll}
 \displaystyle \frac{1}{\nu_2} \int _{\tore}| H_N (\delta \vec{u}) |^2 | \B_2|^2
 &\le \displaystyle \frac{1}{\nu_2} \| H_N (\delta \vec{u}) \|_{L^3}^2 \| \B_2\|_{L^6}^2&\\
 & \displaystyle\le  \frac{1}{\nu_2} \| H_N (\delta \vec{u}) \|_{\vec{V}^{\frac{1}{2}}}^2 \| \B_2\|_{\vec{V}}^2.&
 \end{array}
 \end{align}
\end{subequations}

Hence,
\begin{equation}
 \label{soirabc}
\begin{array}{llll}
 \displaystyle {\frac{d }{2dt}} \int_{\tore} | \delta\vec{u}| ^{2}+{\frac{d }{2dt}} \int_{\tore}  | \delta\B| ^{2} +\displaystyle \frac{\nu_1}{2}   \int_{\tore}| \nabla
\delta\vec{u} |^{2} +\displaystyle \frac{\nu_2}{2} \int_{\tore} | \nabla
\delta\B |^{2} \\
\hskip 1cm \displaystyle \le \frac{1}{\min{(\nu_1,\nu_2)}}\left(
 \| H_N (\delta \vec{u}) \|_{\vec{V}^{\frac{1}{2}}}^2  +\displaystyle   \| H_N (\delta \B) \|_{\vec{V}^{\frac{1}{2}}}^2 \right)\left(\| \vec{u}_2\|_{\vec{V}}^2 +   \| \B_2\|_{\vec{V}}^2    \right)
\\
\hskip 2cm \displaystyle \le \frac{C(\alpha,N)}{\min{(\nu_1,\nu_2)}}\left(\| \vec{u}_2\|_{\vec{V}}^2 +   \| \B_2\|_{\vec{V}}^2    \right) \left(
 \| \delta \vec{u} \|_{\vec{H}}^2  +\displaystyle   \|  \delta \B \|_{\vec{H}}^2 \right).
 \end{array}
\end{equation}

Therefore,
$$
{\frac{d }{2dt}} \left(|| \delta \vec{u} ||_{\vec{H}}^2 + || \delta \B ||_{\vec{H}}^2\right)
 \le C (t) \left(|| \delta\vec{u}||_{\vec{H}}^{2} + || \delta\B||_{\vec{H}}^{2}   \right),
$$
where $\displaystyle C(t) = \frac{C(\alpha,N)}{\min{(\nu_1,\nu_2)}} \left(\| \vec{u}_2\|_{\vec{V}}^2 +   \| \B_2\|_{\vec{V}}^2    \right) \in L^1 ([0,T])$. 
We conclude that $\delta\vec{u} =\delta\B= 0$ thanks to Gr\"{o}nwall's Lemma.\\
In order to deduce the local 
energy equality (\ref{localMHD}), we multiply the equation (\ref{IMHDV1}) with $2\vec{u}\phi,$ the equation  (\ref{IMHDV2}) with $2\B\phi,$
for all $T \in (0,+\infty]$ and for all  non negative fonction  $\phi \in C^{\infty}$ with supp $\phi \subset\subset \tore\times(0,T),$  and then adding them.  The rest can be done in exactly
way as in Lemma \ref{lemmelocalen}, so we omit the details.\\

With smooth initial data we may also prove the following theorem 
\begin{Theorem}
\label{existensssIMHDV1/2}
 For $\theta=1/4$.  Assume ${\vec{u}}^{}_{0} \in \vec{V}^{\frac{1}{2}} $ and $\B_{0} \in \vec{V}^{\frac{1}{2}}
$.  Then for any $T>0$, (\ref{IMHDV}) has a unique regular    solution $ (\vec{u}^{ },\B,{p}^{ }):=(\vec{u}_{\alpha,N },\B_{\alpha,N } ,{p}_{\alpha,N })$, 
 such that 
 $\displaystyle {\vec{u}^{ }}, {\B^{}} \in C([0,T],\vec{V}^{\frac{1}{2}}) \cap L^{2}([0,T],\vec{V}^{\frac{3}{2}}),$
  and 
 $\displaystyle p \in L^{2}([0,T],{H}^{\frac{1}{2}}(\tore)).$ 
\end{Theorem}
\textbf{Proof.}
We only sketch  the proof since  is similar to the one for  Navier-Stokes equations.
The proof  is obtained by taking the inner product of (\ref{IMHDV1}) with $|\nabla|\vec{u},$ (\ref{IMHDV2}) with $|\nabla|B$ and then adding them, the existence of a
solution to Problem $(\ref{IMHDV})$ can be derived thanks to the Galerkin
method. Notice that  ${\vec{u}^{}}, {\B^{}}$    satisfy the following estimates
\begin{equation}
\begin{array}{llllllll}
 \hskip -1cm\displaystyle   \frac{1}{2} \frac{d}{dt}\left(
\|{\vec{u}}^{ }\|_{\vec{V}^{\frac{1}{2}}}^{2}  +
\|{\B}^{}\|_{\vec{V}^{\frac{1}{2}}}^{2}\right) +\min{(\nu_1,\nu_2)} \left(\|\nabla
{\vec{u}}^{ }\|_{\vec{V}^{\frac{1}{2}}}^{2} +\|\nabla
{\B}^{ }\|_{\vec{V}^{{1}{2}}}^{2}\right) \\
\le  \left| \displaystyle  \int_{\tore} H_{N}({\vec{u}})  \nabla \vec{u}^{ }| \nabla |\vec{u}^{ }d\vec{x} \right|+ \left| \displaystyle  \int_{\tore} H_{N}({{\B}})  \nabla {\B^{}}| \nabla |\vec{u}^{} d\vec{x} \right| \\
\hskip 1cm + \left| \displaystyle  \int_{\tore} H_{N}({\vec{u}})  \nabla \B^{\alpha,N }| \nabla |\B^{\alpha,N } d\vec{x} \right|+ \left| \displaystyle  \int_{\tore} H_{N}(\B)  \nabla \vec{u}^{}| \nabla |\B^{} d\vec{x} \right|.
\end{array}
\end{equation}
The first non linear term is estimated by

\begin{equation}
\begin{array}{llllll}
 \left| \displaystyle  \int_{\tore} H_{N}({\vec{u}})  \nabla \vec{u}^{}| \nabla |\vec{u}^{}d\vec{x} \right|
 &\le & \|  H_{N}({\vec{u}})^{}\nabla\vec{u}_{}^{}\|_{\vec{V}^{-{\frac{1}{2}}}} \|\nabla \vec{u}_{}^{} \|_{\vec{V}^{{1}{2}}} \\ 
  &\le& \|  H_{N}({\vec{u}})^{}\nabla\vec{u}\|_{{L}^{\frac{3}{2}}} \|\nabla \vec{u}\|_{\vec{V}^{\frac{1}{2}}} \\ 
   &\le& \|  H_{N}({\vec{u}})^{}\|_{{L}^{6}}\|\nabla\vec{u}\|_{L^2} \|\nabla \vec{u} \|_{\vec{V}^{\frac{1}{2}}} \\
    &\le& \|  H_{N}({\vec{u}})\|_{\vec{V}}\|\nabla\vec{u}\|_{\vec{H}} \|\nabla \vec{u} \|_{\vec{V}^{\frac{1}{2}}}\\
     &\le& \|  H_{N}({\vec{u}})\|_{\vec{V}}\|\vec{u}\|_{\vec{V}^{\frac{1}{2}}}^{\frac{1}{2}} \|\nabla \vec{u} \|_{\vec{V}^{\frac{1}{2}}}^{\frac{3}{2}}.
  \end{array}
\end{equation}
Where we haved used the Sobolev embedding in Theorem \ref{soblevin} together with H\"{o}lder estimate
and  the interpolation inequality between $\vec{V}^{\frac{1}{2}}$ and $\vec{V}^{\frac{3}{2}}$.\\
Similarly, we can estimat the third term by 
\begin{equation}
\begin{array}{llllll}
 \left| \displaystyle  \int_{\tore} H_{N}({\vec{u}})  \nabla \B^{}| \nabla |\B^{ }d\vec{x} \right|
     &\le& \|  H_{N}({\vec{u}})^{}\|_{}\|\B_{}^{}\|_{\vec{V}^{\frac{1}{2}}}^{\frac{1}{2}} \|\nabla \B_{}^{} \|_{\vec{V}^{\frac{1}{2}}}^{\frac{3}{2}}.
  \end{array}
\end{equation}
The second non-linear term is estimated by

\begin{equation}
\begin{array}{llllll}
   \left| \displaystyle  \int_{\tore} H_{N}({{\B}})  \nabla {\B^{}}| \nabla |\vec{u}^{} d\vec{x} \right|
 &\le & \|  H_{N}({\B})^{}\nabla\B_{}^{}\|_{\vec{V}^{-{\frac{1}{2}}}} \|\nabla \vec{u}_{}^{} \|_{\vec{V}^{\frac{1}{2}}} \\ 
  &\le& \|  H_{N}({\B})^{}\nabla\B_{}^{}\|_{{L}^{\frac{3}{2}}} \|\nabla \vec{u}_{}^{} \|_{\vec{V}^{\frac{1}{2}}} \\ 
   &\le& \|  H_{N}({\B})^{}\|_{{L}^{6}}\|\nabla\B_{}^{}\|_{L^2} \|\nabla \vec{u}_{}^{} \|_{\vec{V}^{\frac{1}{2}}} \\
    &\le& \|  H_{N}({\B})^{}\|_{\vec{V}}\|\nabla\B_{}^{}\|_{\vec{H}} \|\nabla \vec{u}_{}^{} \|_{\vec{V}^{\frac{1}{2}}}\\
     &\le& C(\alpha,N)\|  {\B}^{}\|_{\vec{V}^{\frac{1}{2}}}\|\B_{}^{}\|_{\vec{V}}^{} \|\nabla \vec{u}_{}^{} \|_{\vec{V}^{\frac{1}{2}}}^{}.
  \end{array}
\end{equation}

By the same way we obtain 
\begin{equation}
\begin{array}{llllll}
  \left| \displaystyle  \int_{\tore} H_{N}({{\B}})  \nabla {\vec{u}^{}}| \nabla |\B^{} d\vec{x} \right|
     &\le& C(\alpha,N)\|  {\B}^{}\|_{\vec{V}^{\frac{1}{2}}}\|\vec{u}^{}\|_{\vec{V}}^{} \|\nabla \vec{u}_{}^{} \|_{\vec{V}^{\frac{1}{2}}}^{}.
  \end{array}
\end{equation}

Therefore,

\begin{equation}
\begin{array}{llllllll}
\hskip -1cm  \displaystyle   \frac{1}{2} \frac{d}{dt}\left(
\|{\vec{u}}^{ }\|_{\vec{V}^{\frac{1}{2}}}^{2}  +
\|{\B}^{}\|_{\vec{V}^{\frac{1}{2}}}^{2}\right) +\min{(\nu_1,\nu_2)} \left(\|\nabla
{\vec{u}}^{}\|_{\vec{V}^{\frac{1}{2}}}^{2} +\|\nabla
{\B}^{}\|_{\vec{V}^{\frac{1}{2}}}^{2}\right) \\
 \le \|  H_{N}({\vec{u}})^{}\|_{\vec{V}}\|\vec{u}_{}^{}\|_{\vec{V}^{\frac{1}{2}}}^{\frac{1}{2}} \|\nabla \vec{u}_{}^{} \|_{\vec{V}^{\frac{1}{2}}}^{\frac{3}{2}} +\|  H_{N}({\vec{u}})^{}\|_{\vec{V}}\|\B_{}^{}\|_{\vec{V}^{\frac{1}{2}}}^{\frac{1}{2}} \|\nabla \B_{}^{} \|_{\vec{V}^{\frac{1}{2}}}^{\frac{3}{2}}\\
\hskip 1cm +C(\alpha,N)\|  {\B}^{}\|_{\vec{V}^{\frac{1}{2}}}\|\B_{}^{}\|_{\vec{V}}^{} \|\nabla \vec{u}_{}^{} \|_{\vec{V}^{\frac{1}{2}}}^{}+C(\alpha,N)\|  {\B}^{}\|_{\vec{V}^{\frac{1}{2}}}\|\vec{u}^{}\|_{\vec{V}}^{} \|\nabla \vec{u}_{}^{} \|_{\vec{V}^{\frac{1}{2}}}^{}
\end{array}
\end{equation}

By Young inequality,
\begin{equation}
\begin{array}{llllllll}
 \hskip -1cm \displaystyle    \frac{d}{dt}\left(
\|{\vec{u}}^{}\|_{\vec{V}^{\frac{1}{2}}}^{2}  +
\|{\B}^{}\|_{\vec{V}^{\frac{1}{2}}}^{2}\right) + \min{(\nu_1,\nu_2)}\left(\|\nabla
{\vec{u}}^{}\|_{\vec{V}^{\frac{1}{2}}}^{2} +\|\nabla
{\B}^{}\|_{\vec{V}^{\frac{1}{2}}}^{2}\right) \\
\le \displaystyle \frac{c}{\min{(\nu_1,\nu_2)}} \|  H_{N}({\vec{u}})^{}\|_{}^4\|\vec{u}_{}^{}\|_{{\vec{V}^{\frac{1}{2}}}}^{2}+  \displaystyle \frac{c}{\min{(\nu_1,\nu_2)}}\|  H_{N}({\vec{u}})^{}\|_{\vec{V}}^4\|\B_{}^{}\|_{\vec{V}^{\frac{1}{2}}}^{2}\\
\hskip 1cm \displaystyle +\frac{C(\alpha,N)}{\min{(\nu_1,\nu_2)}}\|  {\B}^{}\|_{\vec{V}^{\frac{1}{2}}}^{2}\|\B_{}^{}\|_{\vec{V}}^{2} \displaystyle +\frac{C(\alpha,N)}{\min{(\nu_1,\nu_2)}}\|  {\B}^{}\|_{\vec{V}^{\frac{1}{2}}}^2 \|\vec{u}^{}\|_{\vec{V}}^{2}.
\end{array}
\end{equation}

Therefore,

\begin{equation}
\begin{array}{llllllll}
 \hskip -1cm \displaystyle    \frac{d}{dt}\left(
\|{\vec{u}}^{ }\|_{\vec{V}^{\frac{1}{2}}}^{2}  +
\|{\B}^{}\|_{\vec{V}^{\frac{1}{2}}}^{2}\right) +\min{(\nu_1,\nu_2)} \left(\|\nabla
{\vec{u}}^{}\|_{\vec{V}^{\frac{1}{2}}}^{2} +\|\nabla
{\B}^{ }\|_{\vec{V}^{\frac{1}{2}}}^{2}\right) \\
\le   \displaystyle \|  {\B}^{}\|_{\vec{V}^{\frac{1}{2}}}^2  \left( \frac{c}{\min{(\nu_1,\nu_2)}}\|  H_{N}({\vec{u}})^{}\|_{\vec{V}}^4
+\displaystyle \frac{C(\alpha,N)}{\min{(\nu_1,\nu_2)}}\|\B_{}^{}\|_{\vec{V}}^{2} +\displaystyle \frac{C(\alpha,N)}{\min{(\nu_1,\nu_2)}} \|\vec{u}^{}\|_{\vec{V}}^{2}\right)\\
\hskip 1cm + \displaystyle \frac{c}{\min{(\nu_1,\nu_2)}} \|  H_{N}({\vec{u}})^{}\|_{\vec{V}}^4 \|\vec{u}_{}^{}\|_{\vec{V}^{\frac{1}{2}}}^{2}\\
\hskip 2cm \le \displaystyle C(t) \left(  \|\vec{u}_{}^{}\|_{\vec{V}^{\frac{1}{2}}}^{2}+   \|  {\B}^{}\|_{\vec{V}^{\frac{1}{2}}}^2  \right)
\end{array}
\end{equation}

where $\displaystyle C(t) =  \left( \frac{c}{\min{(\nu_1,\nu_2)}}\|  H_{N}({\vec{u}})^{}\|_{\vec{V}}^4
+\displaystyle \frac{C(\alpha,N)}{\min{(\nu_1,\nu_2)}}\|\B_{}^{}\|_{\vec{V}}^{2} +\displaystyle \frac{C(\alpha,N)}{\min{(\nu_1,\nu_2)}} \|\vec{u}^{}\|_{\vec{V}}^{2}\right) \in L^1 ([0,T])$. 
We conclude that    ${\vec{u}^{}}, {\B^{}} \in L^{\infty}((0,T],\vec{V}^{\frac{1}{2}})\cap L^{2}([0,T],\vec{V}^{\frac{3}{2}})$  thanks to Gr\"{o}nwall's Lemma.\\

The pressure $p^{ }$   is reconstructed from $\bu, H_N(\bu), \B, H_N(\B)$ ( as we are working with periodic boundary conditions)
and its regularity
results from the fact that  $(H_N(\bu)\cdot\nabla)u, (H_N(\B)\cdot\nabla)\B  \in L^2([0,T],\vec{H}^{-{\frac{1}{2}}})$.\\
It remains to prove the uniqueness. Let $(\vec{u}_{1},\B_1, p_{1})$ and $(\vec{u}_{2},\B_2,
p_{2})$ be two solutions, $\delta\vec{u} = \vec{u}_{2}-\vec{u}_{1}$,$\delta\B = \B_{2}-\B_{1}$, $\delta p =
p_{2}-p_{1}$. Then one has

  \begin{equation}
\label{matin1prime}\begin{array}{ll}
 \displaystyle \partial_{t} \delta\vec{u} + (H_{N} (\vec{u}_{1}) \nabla) \delta\vec{u}-(H_{N} (\B_{1}) \nabla) \delta\B  - \nu
\Delta\delta\vec{u} + \nabla\delta p\\
 \hskip 5cm = -(H_{N} (\delta\vec{u}) \nabla) \vec{u}_{2} +(H_{N} (\delta\B) \nabla) \B_{2},  
\end{array}
\end{equation}

  \begin{equation}
\label{matin2prime}\begin{array}{ll}
 \displaystyle \partial_{t} \delta\B +(H_{N} (\vec{u}_{1}) \nabla) \delta\B -  (H_{N} (\B_{1}) \nabla) \delta\vec{u} - \nu
\Delta\delta\B \\
 \hskip 5cm =  (H_{N} (\delta\B) \nabla) \vec{u}_{2}- (H_{N} (\delta\vec{u}) \nabla) \B_{2},
\end{array}
\end{equation}

and $\delta\vec{u} = 0$, $\delta\B =0$  at initial time. One can take $|\nabla|\delta\vec{u} \in  L_{}^{\infty}((0,T],\vec{V}^{-{\frac{1}{2}}})\cap L^{2}([0,T],\vec{V}^{\frac{1}{2}})$
 as test in (\ref{matin1prime}) and $|\nabla|\delta\B \in  L_{}^{\infty}((0,T],\vec{V}^{-{\frac{1}{2}}})\cap L^{2}([0,T],\vec{V}^{\frac{1}{2}})$
 as test in (\ref{matin2prime}). 

Once we obtain  by a similar way as in Theorem \ref{existensssIMHDV}  and Theorem \ref{1} that 

$$
{\frac{d }{2dt}} \left(\| \delta \vec{u} \|_{\vec{V}^{\frac{1}{2}}}^2 + \| \delta \B \|_{\vec{V}^{\frac{1}{2}}}^2\right)
 \le C (t) \left(\| \delta\vec{u}\|_{\vec{V}^{\frac{1}{2}}} ^{2} + \| \delta\B \|_{\vec{V}^{\frac{1}{2}}} ^{2}   \right),
$$
where $\displaystyle C(t)  \in L^1 ([0,T])$. 
We conclude that $\delta\vec{u} =\delta\B= 0$ thanks to Gr\"{o}nwall's Lemma.\\


Next, we will deduce that the deconvolution regularization   give rise to a suitable weak solution to the MHD equations.
 \begin{Theorem}
Let  $(\vec{u}_{\alpha,N},\B_{\alpha,N}, p_{\alpha,N})$ be the solution of  (\ref{dalpha ns}), Then when  $\alpha \rightarrow 0$ and/or $N \rightarrow \infty$,  $(\vec{u}_{\alpha,N},\B_{\alpha,N},, p_{\alpha, N})$ tends to a  weak solution  $(\vec{u}_{},\B, p_{})$ to the  MHD equations. The convergence to $\vec{u}$ and the convergence to $\B$ are weak in $L^2([0,T];\vec{V})\cap L^{\frac{10}{3}}([0,T];L^{\frac{10}
{3}}(\tore)^3)$ and strong in $L^q([0,T];L^q(\tore)^3)$  for all $ 2 \le q <\frac{10}{3}$ and the convergence to $p$ is strong  in $L^q([0,T];L^q(\tore))$  for all $ 1 < q<\frac{5}{3}$ and weak  in  $L^{\frac{5}{3}}([0,T];L^{\frac{5}{3}}(\tore))$.
Furthermore  the solution verifies in addition 
 \begin{equation}
\label{localMHDzero}
\begin{array}{cccc}
\hskip -3cm \displaystyle 2\nu\int_{0}^{T}\int_{\tore}\left(|\nabla \vec{u}_{}|^{2}+|\nabla \B_{}|^{2} \right)\phi \ d\vec{x}dt \\
\hskip 1cm \le \displaystyle \int_{0}^{T}\int_{\tore}| \vec{u}_{}|^{2}\left(\phi_t + \nu_1 \Delta \phi \right) +| \B_{}|^{2}\left(\phi_t + \nu_2 \Delta \phi \right) \  d\vec{x}dt\\
\hskip 1.5cm +\displaystyle\int_{0}^{T}\int_{\tore}  \left( \left(| \vec{u}_{}|^{2}+| \B_{}|^{2}\right) {\vec{u}_{}}+ 2p\vec{u}_{}\right) \cdot \nabla \phi   \  d\vec{x}dt\\
 \hskip 2cm 
 \displaystyle -\displaystyle 4\int_{0}^{T}\int_{\tore}  \left( \vec{u}_{} \B_{}\right){\B_{}}\cdot \nabla \phi\  d\vec{x}dt 
 \end{array}
\end{equation}
for all $T \in (0,+\infty]$ and for all  non negative fonction  $\phi \in C^{\infty}$ and  supp $\phi \subset\subset \tore\times(0,T).$
\end{Theorem}
\textbf{Proof.}
 As in Lemma \ref{conlp} we can show that  
 \begin{align}
H_N({\bu_{\alpha,N}})  &\rightarrow \bu &&\textrm{strongly in  }
L^p(0,T;L^p(\tore)^3) \textrm{ for all  }   2 \le p < \frac{10}{3}.\\
H_N({\B_{\alpha,N}})  &\rightarrow \B &&\textrm{strongly in  }
L^p(0,T;L^p(\tore)^3) \textrm{ for all  }   2 \le p < \frac{10}{3}.
\end{align}
  The above $L^p$ convergence combined with the fact that  $\vec{u}_{\alpha,N}$ and $\B_{\alpha,N} $  belong to  the enegy space of the solutions of the Navier-Stokes equations and the  Aubin-Lions compactness  Lemma allow us to take the limit $ \alpha \rightarrow 0$ or/and $ N \rightarrow  \infty$ in (\ref{localMHD}). 
 The rest can be done in exactly
way as in \cite{LL08}, so we omit the details.\\

\end{document}